\documentclass{article}       
\usepackage{amsmath,amsthm,amsfonts,psfrag,amssymb,dsfont,mathrsfs,mathtools,comment}
\usepackage{authblk}
\usepackage{enumitem,linegoal}
\usepackage{tikz-cd}
\usepackage{comment}
\usepackage{hyperref}
\hypersetup{
     colorlinks=true,
     linkcolor=red,
     filecolor=blue,
     citecolor = blue,      
     urlcolor=blue,
     }
\usepackage{cleveref}

\theoremstyle{plain}
\newtheorem{theorem}{Theorem}
\crefname{thm}{theorem}{theorems}
\newtheorem{thmintro}{Theorem}

\crefname{thmintro}{theorem}{theorems}
\newtheorem{obs}[theorem]{Observation}
\newtheorem{corintro}[thmintro]{Corollary}

\crefname{corintro}{corollary}{corollaries}
\newtheorem{lemma}[theorem]{Lemma}
\newtheorem{corollary}[theorem]{Corollary}
\newtheorem{proposition}[theorem]{Proposition}

\theoremstyle{definition}
\newtheorem{definition}[theorem]{Definition}
\newtheorem{remark}[theorem]{Remark}

\newtheorem{question}[theorem]{Question}

\def\co{\colon\thinspace}
\def\coeq{\coloneqq\thinspace}

\newcommand{\RR}{\mathbb{R}}

\newcommand{\TT}{\mathbb{T}}
\newcommand{\NN}{\mathbb{N}}
\newcommand{\ZZ}{\mathbb{Z}}
\newcommand{\DD}{\mathbb{D}}
\newcommand{\QQ}{\mathbb{Q}}
\newcommand{\SSS}{\mathbb{S}}

\newcommand{\calM}{\mathcal{M}}

\newcommand{\bigslant}[2]{{\raisebox{.2em}{$#1$}\left/\raisebox{-.2em}{$#2$}\right.}}

\DeclareMathOperator{\image}{image}

\DeclareMathOperator{\Int}{int}

\DeclareMathOperator{\Id}{Id}

\newcommand{\neigh}{\mathcal{N}}

\newcommand{\omegahat}{\widehat{\omega}}

\newcommand{\calB}{\mathcal{B}}

\usepackage{authblk}

\newcommand{\OBD}{OBD}
\newcommand{\BO}{BO}

\newcommand{\lambdastd}{\lambda_{std}}
\newcommand{\omegastd}{\omega_{std}}
\newcommand{\xistd}{\xi_{std}}

\newcommand{\Sigmahat}{\widehat{\Sigma}}

\newcommand{\calI}{\mathcal{I}}
\newcommand{\psihat}{\widehat{\psi}}
\newcommand{\intW}{\mathring{\Sigma}}
\newcommand{\PS}{\mathcal{PS}}
\newcommand{\norm}[1]{\left\vert#1\right\vert}
\newcommand{\Ctop}{C_{top}}
\newcommand{\Cbot}{C_{bot}}

\newcommand{\calMbar}{\overline{\calM}}
\newcommand{\calMstarbar}{\overline{\calM}_*}
\newcommand{\Wcap}{W_{cap}}

\newcommand{\tildej}{\widetilde{j}}

\newcommand{\tildeW}{\widetilde{W}}

\begin{document}

\title{Bourgeois contact structures: \\
tightness, fillability and applications}

\author{Jonathan Bowden\footnote{Monash University, Melbourne, Australia.  Current: Universit\"at Regensburg, Regensburg, Germany. Email: \url{jonathan.bowden@mathematik.uni-regensburg.de}}, 
Fabio Gironella\footnote{Alfred Renyi Institute, Budapest, Hungary. Current: Humboldt University, Berlin, Germany. Email: \url{gironelf@math.hu-berlin.de}}, 
Agustin Moreno\footnote{Universit\"at Augsburg, Augsburg, Germany. Current: Uppsala Universitet, Uppsala, Sweden. Email: \url{agustin.moreno2191@gmail.com}}}

%\affil[1]{Monash University, Melbourne, Australia.  Current: Universit\"at Regensburg, Regensburg, Germany. Email: \url{jonathan.bowden@mathematik.uni-regensburg.de}}
%\affil[2]{Alfred Renyi Institute, Budapest, Hungary. Current: Humboldt University, Berlin, Germany. Email: \url{gironelf@math.hu-berlin.de}}
%\affil[3]{Universit\"at Augsburg, Augsburg, Germany. Current: Uppsala Universitet, Uppsala, Sweden. Email \url{agustin.moreno2191@gmail.com}}

\date{}
\maketitle
\begin{abstract}
	Given a contact structure on a manifold $V$ together with a supporting open book decomposition, Bourgeois gave an explicit construction of a contact structure on $V \times \mathbb{T}^2$. 
	We prove that all such structures are universally tight in dimension $5$, independent on whether the original contact manifold is itself tight or overtwisted.
	
	In arbitrary dimensions, we provide obstructions to the existence of strong symplectic fillings of Bourgeois manifolds. This gives a broad class of new examples of weakly but not strongly fillable contact $5$--manifolds, as well as the first examples of weakly but not strongly fillable contact structures in all odd dimensions. These obstructions are particular instances of more general obstructions for $\SSS^1$-invariant contact manifolds.
	
	We also obtain a classification result in arbitrary dimensions, namely that the unit cotangent bundle of the $n$-torus has a unique symplectically aspherical strong filling up to diffeomorphism.
\end{abstract}

%\normalsize
\tableofcontents

\section{\textbf{Introduction}} 

In \cite{Bo}, Bourgeois showed that, whenever $(V,\xi)$ is a contact manifold endowed with a supporting open book decomposition (which always exist by work of Giroux \cite{Gir}), then the manifold $V \times \mathbb{T}^2$ carries a natural contact structure. 
The main motivation behind such a construction was the problem of the existence of contact structures on higher--dimensional manifolds. 
For instance, it showed that every odd dimensional torus admits contact structures, a problem that has been open since Lutz \cite{Lutz79} proved that $\mathbb{T}^5$ is contact, more than 20 years before. 

It was not until recently that Borman-Eliashberg-Murphy \cite{BEM} proved that contact structures in higher-dimensions actually exist in abundance (i.e.\ whenever the obvious topological obstructions disappear) by generalizing Eliashberg's \cite{Eli89} notion of \emph{overtwistedness}, as well as the $h$--principle that comes with it, to higher dimensions. Overtwisted contact manifolds are topological/flexible in nature, and most of the associated contact-topological invariants (e.g.\ those coming from holomorphic curves) simply vanish. As a result, it has become relevant to find examples of high--dimensional contact structures beyond the overtwisted ones, which are more geometric/rigid and which potentially have rich associated invariants. Contact structures which are not overtwisted are usually referred to as \emph{tight}.

The construction in \cite{Bo} actually fits very well in this setting, as it is both very explicit and yields a very broad class of contact manifolds in arbitrary odd dimensions, with remarkable properties. For instance, \cite{Pre07} used it to construct the first examples of high dimensional (closed) contact manifolds admitting a Plastikstufe, as defined in \cite{Nie06}, which is equivalent to overtwistedness \cite{CMP,Hua}. In \cite{BCS14}, the authors also used it to construct contact structures on the product of a contact manifold with the $2$--sphere. More recently, Lisi, Marinkovi\'{c} and Niederkr\"uger \cite{LMN} started the systematic study of the Bourgeois construction, in particular studying its fillability properties, and in this paper we continue this line of research. 

In what follows we will use the following notation. Given an abstract open book $(\Sigma^{2n},\phi)$ and the associated contact $(2n+1)$--manifold $OBD(\Sigma,\phi)$, we denote by  $BO(\Sigma,\phi)$ the contact manifold obtained via the Bourgeois construction \cite{Bo}. Smoothly, $BO(\Sigma,\phi)=OBD(\Sigma,\phi)\times \TT^2$ and we refer to \Cref{sec:bourgeois_constr} for further details.

\medskip

\textbf{Tightness in dimension $5$.} We begin by addressing the natural question of whether a given Bourgeois contact structure is tight or overtwisted. 

In \cite{LMN}, the authors give examples, in every odd dimension, of an overtwisted $(V,\xi)$ such that the associated Bourgeois contact manifold is tight. 
Moreover, in \cite{Gi}, it was shown that if $V$ is a $3$-manifold with non-zero first Betti number, then there exists a supporting open book such that the associated Bourgeois contact structure is (hyper)tight.
In this paper, we prove that, at least in dimension $5$, these are particular instances of a more general fact. 
Namely, $5$-dimensional Bourgeois contact structures are rigid, inherently geometric objects, independently of the rigid or flexible nature of $(V,\xi)$:
\begin{thmintro}[Tightness]
	\label{thm:tight_bourgeois}
	For every abstract open book $(\Sigma^2,\phi)$, the contact $5$--manifold $BO(\Sigma,\phi)$ is (universally) tight.
\end{thmintro}
  
Recall that universally tight means that the universal cover is tight. 
In particular, universal tightness implies tightness. The fact that $5$--dimensional Bourgeois contact structures are universally tight is a simple consequence of the fact that they are tight and that finite covers on either factor of the product again yield Bourgeois contact structures.

While there are many ways of making contact/symplectic manifolds more flexible (e.g.\ in dimension $3$, by adding a Lutz twist; in any odd dimension, by taking the connected sum with an overtwisted sphere; or in dimension at least $6$, by taking the flexibilization of a Weinstein manifold), \Cref{thm:tight_bourgeois} says that the Bourgeois construction can be interpreted as a ``tightening'' procedure. To our knowledge, there is currently no other procedure with an analogous property.

Moreover, the above result is sharp with respect to taking branched covers. 
Namely, if $V$ is an overtwisted contact $3$-manifold (or, more generally, $(2n+1)-$manifold), then branched covers $V\times \Sigma_g$ of $V\times \mathbb{T}^2$ are overtwisted; here, $\Sigma_g$ is the orientable surface of genus $g\geq 1$, seen as a degree $g$ branched cover of $\TT^2$ over two points, and the branched cover is obtained by product with the identity on $V$. 
The fact that this is true for $g$ big enough follows from the results in \cite{NiePre10};
that this holds already from $g=2$ is moreover a consequence of an argument of Massot and Niederkr\"{u}ger, based on ideas from \cite{Pre07} (cf.\ \cite[Theorem I.5.1]{NieNotes};
the interested reader can also consult \cite[Observation 5.9]{Gi} for details.)

\medskip

\textbf{Symplectic fillability.} Another important problem in contact topology that is related to the flexible/rigid classification of contact structures, is that of characterizing which contact structures admit symplectic fillings. 
Indeed, one of the first results of the theory of holomorphic curves is that symplectically fillable contact manifolds are tight.

Given \Cref{thm:tight_bourgeois}, it is natural to wonder whether Bourgeois contact structures are (at least weakly) symplectically fillable. 
While there exist partial results in this direction \cite{MNW,LMN}, a complete answer is yet to be found.
%See also \Cref{subtle} for subtleties in finding weak fillings.

In order to find obstructions to strong fillings, we shall prove a more general statement concerning $\SSS^1$-invariant contact structures on manifolds of the form $V^{2n} \times \SSS^1$. Such a contact structure is determined by a splitting of the base $V$ into \emph{ideal} Liouville domains $V=V_+ \cup \overline{V}_-$ glued along a contact manifold $N:=\partial V_+ = \partial V_-$ (see \cite{GirouxIdealLiouvDom,MNW} or Section \ref{sec:general_obstr_fillab} below). 
The subsets $V_\pm$ correspond to the positive resp.\ negative regions when $V$ is viewed as a convex hypersurface, and $N$ to its dividing set. With this notation we then have the following homological and homotopical criteria in the case that the contact structure is strongly fillable.
\begin{thmintro}[Fillability of $\SSS^1$-invariant contact structures]\label{thm:homology_injection_boundary}
Suppose that $V^{2n} \times \SSS^1$ admits a strongly fillable $\SSS^1$-invariant contact structure with induced convex splitting $V =V_+ \cup \overline{V}_-$ and let $N:=\partial V_+ = \partial V_-$ be the dividing set, which we assume to be connected.
Let also $(W^{2n+2},\omega)$ be a strong filling of $V\times \SSS^1$ and consider the natural inclusions and induced maps on homology:
\begin{equation*}
%\label{eqn:thm_homology_injection_page_general_case}
N \hookrightarrow  V_\pm \hookrightarrow W \textrm{ and } H_*(N,\mathbb{Q}) \stackrel{I_\pm}{\longrightarrow}  H_*( V_\pm,\mathbb{Q}) \longrightarrow H_*( W,\mathbb{Q}) .
\end{equation*}
Then the second inclusion induces an injection on the image $Im\left( I_\pm \right)$ in rational homology.
In particular, if $V_\pm$ are Weinstein, then the inclusion $N \hookrightarrow W$ induces an injection in rational homology in all degrees strictly less than $n$. 

In the case where $W$ is semi-positive, we also have a surjection
$$
\pi_1(V\times \mathbb{S}^1)\twoheadrightarrow \pi_1(W)
$$
of fundamental groups.
\end{thmintro}

%The above result may be seen as a version of the Giroux torsion obstruction introduced by Massot-Niederkr\"uger-Wendl \cite{MNW}, which can be recast as the statement that if the number of connected components of $N$ is at least two, then $V\times\mathbb{S}^1$ does not admit a (connected) filling $W$. This corresponds to the degree zero case of the statement in \Cref{thm:homology_injection_boundary}, i.e.\ that concerning $H_0(N,\mathbb{Q})$. Note however that this theorem also provides obstructions if $N$ is assumed connected.
%For simplicity, we assume above that the dividing set is connected, as this also highlights the distinction from the case considered in \cite{MNW}.

Since Bourgeois contact structures are $\TT^2$-invariant, they are in particular $\SSS^1$-invariant and one obtains a convex splitting with pieces of the form $V_\pm = \Sigma \times D^* \SSS^1$. 
Then %it turns out that strong fillability for Bourgeois manifolds is quite a restrictive condition. Indeed, 
applying \Cref{thm:homology_injection_boundary}, we obtain the following result, which imposes strong topological restrictions on the filling as well as the original Bourgeois contact manifold:

\begin{thmintro}[Fillability of Bourgeois contact structures]\label{thm:homology_injection_page}
Suppose that a $(2n+1)$-dimensional Bourgeois contact manifold $BO(\Sigma,\phi)$ is strongly symplectically filled by $(W,\omega)$. 
Then the natural inclusion of a page into the filling, given by the composition of inclusions 
\begin{equation}
\label{eqn:thm_homology_injection_page}
\Sigma \hookrightarrow OBD(\Sigma,\phi)\times\{pt\} \hookrightarrow BO(\Sigma,\phi)=\partial W \hookrightarrow W,
\end{equation}
%with $pt\in\TT^2$, 
induces an injection in rational homology. 

Moreover, the composition of the natural inclusions $\TT^2\hookrightarrow BO(\Sigma,\phi)\hookrightarrow W$ also induces an injection in rational homology. 
Lastly, if $W$ is semi-positive, then we have a surjection $ \pi_1(\Sigma\times \TT^2)\twoheadrightarrow \pi_1(W)$.
\end{thmintro}

\noindent There are many topological situations where the above criteria can be immediately applied to obstruct strong fillability of certain Bourgeois contact structures, several of which are described below. We also have:
\begin{remark}(Monodromy restrictions).
\label{rk:phi_is_id_homology}
The fact that $\Sigma\hookrightarrow W$ induces an injection in rational homology implies that the same is true for the natural inclusion $\Sigma\hookrightarrow OBD(\Sigma,\phi)$, as $\Sigma\hookrightarrow W$ can be written as the composition in \Cref{eqn:thm_homology_injection_page}.
What is more, this also implies that $\phi_*=\Id$ on $H_*(\Sigma;\mathbb Q)$, i.e.\ $\phi$ lies in the Torelli group: indeed, as $\Sigma\hookrightarrow OBD(\Sigma,\phi)$ is the composition $\Sigma \hookrightarrow \Sigma_\phi \hookrightarrow OBD(\Sigma,\phi)$ where $\Sigma_\phi$ is the mapping torus part of the open book, one can appeal to the long exact sequence of a mapping torus \cite[Example 2.48]{HatAlgTop} to deduce that $\Sigma\hookrightarrow \Sigma_\phi$ is injective in rational homology if and only if $\phi_*=\Id$ on $H_*(\Sigma;\mathbb{Q})$.
\end{remark}

\begin{remark}
To deal with strong fillings in general, we need to use the polyfold machinery of Hofer--Wysocki--Zehnder \cite{HWZ_GW} (for the case of \emph{closed} spheres). 
However, this is not necessary if further technical assumptions on the filling are imposed, e.g.\ asphericity or semi-positivity.
%(In the proof in \Cref{sec:general_obstr_fillab} we will treat the semi-positive case separately).
\end{remark}

\noindent \textbf{Homology of aspherical fillings.} In the case where a Bourgeois manifold is assumed to admit a \emph{symplectically aspherical} filling, we then obtain the following significant strengthening of  \Cref{thm:homology_injection_page}, which determines the homology of any such filling:
\begin{thmintro}[Aspherical case]\label{thm:hom_apsherical_filling_Bour}
Suppose that a $(2n+1)$-dimensional Bourgeois contact manifold $BO(\Sigma,\phi)$ admits a symplectically aspherical filling $(W,\omega)$. 
Then the natural inclusion $\Sigma \times \TT^2\hookrightarrow W$ into the filling induces an isomorphism in integral homology.
\end{thmintro}

If $\phi$ is symplectically isotopic to $\Id$, $OBD(\Sigma,\phi)$ admits a subcritical Stein filling. %\cite{Ci02}. 
 In this case the associated Bourgeois contact manifold admits a Stein filling according to \cite[Theorem A.b]{LMN}. 
Moreover, this filling is smoothly of the form $X\times \mathbb{T}^2$, where $X=\Sigma\times \mathbb{D}^2$ is the subcritical filling of the original manifold.
\Cref{thm:homology_injection_page}, as well as \Cref{rk:phi_is_id_homology} and \Cref{thm:hom_apsherical_filling_Bour}, then suggest that perhaps it is always the case that strong fillability for Bourgeois manifolds implies $\phi=\Id$ (at least smoothly), and that the filling is (again, at least smoothly) the standard one, in arbitrary dimensions. 
The fact that the filling ``remembers'' both the $\mathbb{T}^2$-factor and the page $\Sigma$ (at least homologically), as well as the fact that, in the aspherical case, the homology of the filling is the expected one are evidence in this direction. 
In particular, strong fillability for Bourgeois manifolds might be equivalent to Stein fillability; this is indeed true in many cases, as discussed below.

\medskip

\textbf{Fillability in dimension $5$.} Examples of weakly but not strongly fillable contact structures in dimension $3$ are well-known. 
The first examples of such contact manifolds in higher dimensions were obtained in dimension $5$ on manifolds also diffeomorphic to a product of a $3$-manifold with a $2$-torus cf.\ \cite[Theorem E]{MNW}. 
These examples are associated to contact $3$-manifolds arising from so-called \emph{Liouville pairs} %(see \cite[Definition 1]{MNW})
%These are in particular co--fillable, hence necessarily \emph{not} planar by a result of Etnyre \cite{Etn04} (i.e.\ the page of any supporting open book cannot have genus zero). 
and and are known to exist only on very specific $3$-manifolds (see \cite{MNW} and references therein). %not including $\SSS^3$ (cf.\ \cite[Remark 8.7]{MNW}).
Below, we provide a large class of new $5$-dimensional examples that arise via Bourgeois contact structures. %including examples arising from planar contact manifolds (\Cref{cor:pos_dehn_twist_no_asph_fill} below).

Indeed, as a first immediate corollary of \Cref{thm:homology_injection_page} in dimension $5$, we obtain:

\begin{corintro}[Rational homology $3$-spheres]
\label{cor:stein_strong_equiv_homol_spheres}
Suppose that $V=OBD(\Sigma,\phi)$ is a $3$-dimensional rational homology sphere, i.e.\ $H_1(V;\mathbb{Q})=0$. 
Then the Bourgeois contact manifold $BO(\Sigma,\phi)$ is strongly symplectically fillable if and only if the page $\Sigma = \DD^2$ is a disc (and $\phi=\Id$), in which case it is actually Stein fillable and $V$ is $\mathbb{S}^3$. 
\end{corintro}

This then yields many examples of weakly but not strongly fillable contact manifolds in dimension $5$. For example Legendrian surgery on any smoothly non-trivial Legendrian knot in $\SSS^3$ gives a Stein fillable contact structure on a rational homology sphere, smoothly different from the standard $3$-sphere. The corresponding Bourgeois contact structure is weakly fillable by \cite[Theorem A.a]{LMN} (cf.\ \cite[Example 1.1]{MNW}).

\medskip

We now consider the \emph{planar} case, i.e.\ the case where the page of the open book of the original $3$-manifold has genus zero. 
In the $3$-dimensional situation, strong symplectic fillings of contact structures supported by planar open books exist in abundance, and in fact are in $1$-$1$ correspondence with the factorizations of the monodromy into products of positive Dehn-twists \cite{Wen10b}. Namely, contact $3$--manifolds admitting a supporting open book with planar page and monodromy a product of positive Dehn twists are (precisely, by \cite{Wen10b}) the convex boundaries of symplectic Lefschetz fibrations. In particular, in the planar case, strong fillability and Stein fillability are equivalent. However, if we apply the Bourgeois construction to a planar contact $3$-manifold, and consider fillings of the resulting contact $5$-manifold, the situation is surprisingly more rigid:

\begin{thmintro}[Planar case]
\label{cor:pos_dehn_twist_no_asph_fill} 
	Let $(\Sigma^2,\phi)$ be an abstract open book with $\Sigma$ of genus zero and $\phi$ a non-trivial product of Dehn twists, all of the same sign.
	Then, the contact $5$-manifold $BO(\Sigma,\phi)$ is weakly but not strongly fillable.
\end{thmintro} 

The fact that the examples in \Cref{cor:pos_dehn_twist_no_asph_fill} are weakly fillable follows again from \cite[Theorem A.a]{LMN} or \cite[Example 1.1]{MNW}; and the fact that the conclusion also holds for products of all \emph{negative} Dehn twists, from \cite[Theorem  B]{LMN}. %(cf.\ \Cref{subtle} below).
In the particular case of the annulus $\Sigma=D^*\SSS^1$, whose mapping class group is generated by the Dehn twist $\tau$ along the zero section, we get $BO(D^*\SSS^1,\tau^k)$ is strongly fillable if and only if $k=0$; see \Cref{thm:BOfill} below for a higher-dimensional version.

Note that, in the planar case, the monodromy $\phi$ is necessarily the identity in homology of the page, since this is generated by the boundary loops, along which $\phi$ is trivial by assumption. 
In particular, the condition on $\phi$ given by \Cref{rk:phi_is_id_homology} is not restrictive. However, one can prove the following, by applying \Cref{thm:homology_injection_page}:

\begin{corintro}
\label{cor:commutator}
If $\Sigma$ is planar and $BO(\Sigma,\phi)$ is strongly fillable, then $\phi$ lies in the commutator subgroup of the mapping class group (rel.\ boundary).   
\end{corintro} 

This implies, for instance, that whenever $\Sigma$ is a pair of pants (whose mapping class group is abelian) and $BO(\Sigma,\phi)$ is strongly fillable, then $\phi=\Id$. Note that whilst mapping class groups of higher genus surfaces are {\em perfect}, meaning that any element can be written as a product of commutators, this is not the case in the planar case. In fact, any non-trivial product of positive Dehn twists will \emph{not} lie in the commutator subgroup. This is because any product of positive Dehn twists gives a non-trivial positive braid after identifying all but one appropriately chosen boundary component to (marked) points, and positive braids survive in the abelianization of the mapping class group of the marked disk, which is infinite cyclic (cf.\ \cite[p.\ 252]{FM}). In particular, one can then deduce \Cref{cor:pos_dehn_twist_no_asph_fill} from \Cref{cor:commutator} via this observation.

\medskip

\medskip

\textbf{Fillability in higher dimensions.} As a further consequence of \Cref{thm:homology_injection_page} we also obtain the following, which gives a plethora of weakly but not strongly fillable contact structures in arbitrary dimension. 
\begin{corintro}[Stabilizations]
\label{cor:high_dim_stab_not_fill}
 Let $OBD(\Sigma,\phi)$ be a $(2n-1)-$dimensional contact manifold, and $OBD(\Sigma_+,\phi_+)$ obtained by a single positive stabilization. Then the Bourgeois manifold $BO(\Sigma_{+},\phi_{+})$ is not strongly symplectically fillable. In particular, if $OBD(\Sigma,\phi)$ is weakly fillable, then $BO(\Sigma_{+},\phi_{+})$ is weakly but not strongly fillable.
\end{corintro}
The fact that such manifolds are not subcritically Stein fillable was already observed in  \cite[Corollary 1.4]{LMN} and hence the above can be viewed as a significant strengthening of this. Note that in \Cref{cor:high_dim_stab_not_fill} we use that the Bourgeois construction preserves weak fillability. As a consequence we obtain the first known examples of weakly but not strongly fillable contact structures in all dimensions.

See \Cref{thm:bou_tau_high_dim} and \Cref{thm:BOfill} below for more concrete examples. 

\medskip

\textbf{Tight but non-fillable contact $5$-folds.} Recall that an almost contact structure on $V^{2n+1}$ is a hyperplane field $\xi$ equipped with a complex structure $J\co \xi\to \xi$. 
On a $3$-manifold, this simply reduces to the data of a homotopy class of oriented $2$-plane fields. Using Eliashberg's classification of overtwisted contact structures in dimension $3$ \cite{Eli90}, we can represent any almost contact structure on a $3$-manifold by an overtwisted contact structure, which is supported by an open book by Giroux \cite{Gir}. If we perform a single \emph{positive} stabilization, the contact structure is unchanged up to isotopy, and hence the almost contact structure is also unchanged up to homotopy. 
Applying \Cref{cor:high_dim_stab_not_fill} to the stabilized open book, and combining with \Cref{thm:tight_bourgeois}, we then conclude:

\begin{corintro}\label{cor:almost}
	For an almost contact structure $(M^3,\eta)$ on a closed $3$-manifold, there exists a universally tight but not strongly fillable contact structure on $M\times \mathbb{T}^2$ which is homotopic to the product almost contact structure $\eta\oplus T\mathbb{T}^2$.
\end{corintro}

We remark that the contact structure in \Cref{cor:almost} can sometimes be weakly fillable (e.g.\ in the case where the monodromy of the planar open book is a product of all \emph{negative} Dehn twists, for which we do not need to positively stabilize; cf.\ \Cref{cor:pos_dehn_twist_no_asph_fill}).

\medskip

\textbf{Special higher-dimensional examples.} We also consider fillings for certain specific examples of higher dimensional Bourgeois contact manifolds. 

\medskip

\textbf{The family $BO(D^*\SSS^n,\tau^k)$. }%, $n\geq 1$, $k\in \mathbb{Z}$.
 We first consider the case of $BO(D^*\SSS^n,\tau)$, where $\tau$ is the Dehn--Seidel twist on $D^*\SSS^n$, and as a special case of \Cref{cor:high_dim_stab_not_fill}, or alternatively as a direct consequence of \Cref{thm:homology_injection_page}, we give a negative answer to \cite[Question 1.6]{LMN}:
\begin{thmintro}
	\label{thm:bou_tau_high_dim}
	The Bourgeois contact manifold $BO(D^*\SSS^n,\tau)$ is weakly but not strongly symplectically fillable. 
\end{thmintro}

\noindent More generally, from \Cref{thm:homology_injection_page} one easily obtains that the Bourgeois contact manifold $BO(D^*\SSS^n,\tau^k)$ is not symplectically fillable for most values of $k$ and $n$:

\begin{thmintro}
\label{thm:BOfill} 
For $n\geq 1$, define the subset
$$
BOFill(n)=\left\{k \in \mathbb{Z}: BO(D^*\SSS^n,\mathbb{\tau}^k) \mbox{ is strongly fillable}\right\}.
$$
Then $BOFill(n)$ is a subgroup of $\ZZ$. 
Denoting a generator by $k_0(n)$, we moreover have the following: 
if $n$ is odd, then $k_0(n)=0$, i.e.\ $BOFill(n)$ is the trivial group, and
if $n$ is even, then $k_0(n)$ is even.
\end{thmintro}

\begin{remark}
The fact that $k_0(n)$ is even if $n$ is also is related to the fact that, in this case, the Dehn-Seidel twist has finite order as a smooth map; see \cite{Keating2021OnTO} for the precise orders. We remark that it is conceivable that $BOFill(n)$ is always the trivial group, which would be a stronger result, although the group itself would no longer be interesting.
\end{remark}

On the other hand, the contact manifolds $BO(D^*\SSS^n,\tau^k)$ admit weak fillings for \emph{every} $k\in \mathbb{Z}$ \cite[Theorem A.a]{LMN}. We then obtain infinitely many weakly but not strongly fillable examples from this family alone; cf.\ \Cref{thm:bou_tau_high_dim}.

\medskip

\textbf{Unit cotangent bundle of $\mathbb{T}^n$.} We determine the diffeomorphism type of the strong symplectically aspherical fillings of the unit cotangent bundle $S^*\mathbb{T}^n$ of $\mathbb{T}^n$ for $n\geq 2$ with its standard contact structure $\xi_{std}$, induced by the restriction of the standard Liouville form $\lambda_{std}$ on the unit disc cotangent bundle $D^*\mathbb{T}^n$ to its boundary $S^*\mathbb{T}^n$. 

In fact, $(S^*\mathbb{T}^n,\xi_{std})$ is none other than $BO(D^*\TT^{n-2},\Id)$. 
This follows from the following. The contact manifold  $OBD(D^*\mathbb{T}^{n-2},\Id)$ is the convex boundary of the subcritical Stein manifold 
\begin{equation*}
W=\underbrace{D^*\SSS^1\times\dots\times D^*\SSS^1}_{n-2 \text{ times}}\times \mathbb{D}^2 \text{ .}
\end{equation*}
According to \cite[Theorem A.b]{LMN}, $BO(T^*\TT^{n-2},\Id)$ is then the convex boundary of the Stein manifold $\left(\prod_{i=1}^{n-2}D^*\mathbb{S}^1\right) \times D^*\TT^2=D^*\TT^n$, with its split Stein structure, which is just the standard one. 

We then prove the following uniqueness result, which is a strengthening of \Cref{thm:hom_apsherical_filling_Bour} for the case of $S^*\mathbb{T}^n$:

\begin{thmintro}
	\label{thm:application}
	The contact manifold $(S^*\mathbb{T}^n,\xi_{std})$, $n\geq 3$, has a unique strong symplectically aspherical filling up to diffeomorphism.
\end{thmintro}

\Cref{thm:application} is a smooth higher--dimensional version of a result by Wendl \cite[Theorem 4]{Wen10b} who classified symplectic fillings of $S^*\mathbb{T}^2$ up to symplectic deformation, which in turn generalized a weaker homeomorphism classification of Stipsicz \cite[Theorem 1.6]{Stipsicz}.
In dimension $5$ (i.e.\ $n=3$), \Cref{cor:pos_dehn_twist_no_asph_fill} and \Cref{thm:application} then give a complete \emph{smooth} characterization of symplectically aspherical strong fillings for the Bourgeois contact $5$--manifolds associated to any open book with page $D^*\SSS^1$; i.e.\ the only strongly fillable case is the trivial monodromy case, in which case the filling is smoothly unique. Note that, according to \cite{LMN}, all examples with page $D^*\SSS^1$ are weakly fillable.

\begin{remark} 
	\Cref{thm:application} has also been independently obtained by Geiges--Kwon--Zehmisch \cite{GKZ}. 
	While our original proof made use of punctured holomorphic curves, the current one uses now closed holomorphic spheres, as kindly suggested to us by a referee in order to simplify the arguments. 
	Hence it now follows a similar line of argument as the proof in \cite{GKZ}.
\end{remark}

\textbf{Outline of the proofs.} For convenience of the reader, we outline the main arguments of the proofs of \Cref{thm:tight_bourgeois,thm:homology_injection_page,thm:application}.

\medskip
\textsc{Tightness in dimension $5$.} 
The proof of \Cref{thm:tight_bourgeois} involves some geometric group theory and hyperbolic geometry as well as some holomorphic curve techniques. 

\begin{figure}[t]
	\centering
	\includegraphics[width=0.35\linewidth]{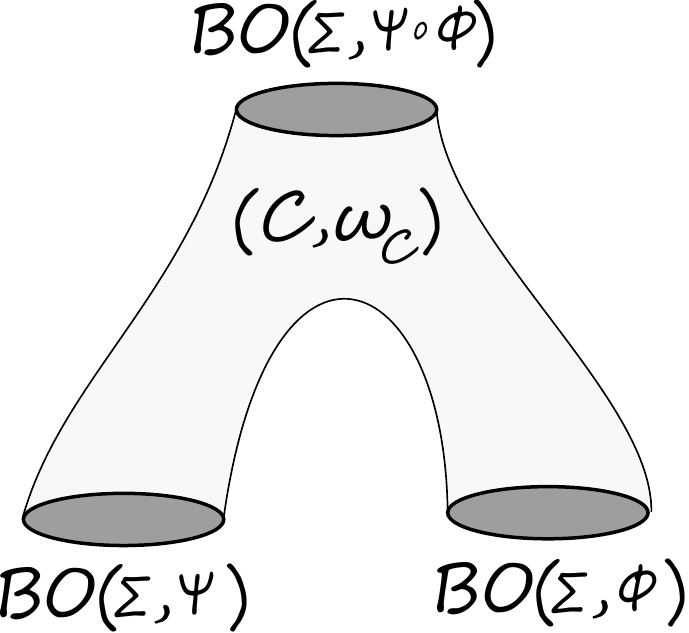}
	\caption{\normalsize The pseudo-Liouville cobordism $(C,\omega_C)$ given in \Cref{thm:strong_cobord_bourgeois}.}
	\label{StabKlfig} 
\end{figure}

The first ingredient is the construction of a strong symplectic cobordism between Bourgeois contact structures; this is done in \Cref{sec:stab_klukas}. 
More precisely, \Cref{thm:strong_cobord_bourgeois} is a ``stabilized'' version of the analogous result for open books, which was proven (independently) in \cite{Advek,Klu12}; see \Cref{StabKlfig}.
We point out that, while the symplectic form on the strong cobordism of \Cref{thm:strong_cobord_bourgeois} is exact, the Liouville vector field associated to the global primitive is \emph{not} inwards pointing along the negative ends. 
We shall refer to a strong symplectic cobordism with an exact symplectic form as \emph{pseudo-Liouville}.

For ``most'' cases of surfaces (the rest are dealt with case by case), standard results from low-dimensional topology then allow to write any monodromy as a composition such that the contact structures on the negative ends of the cobordism in \Cref{StabKlfig} are hypertight; see \Cref{hypertightfactorization}. 
Then, a standard application of the holomorphic curve machinery \`a-la \cite{AH09,Hof93,Nie06} gives a holomorphic plane in the symplectization of one of the negative ends starting with a Bishop family associated to a Plastikstufe in the positive end. While bubbles are ruled out by exactness, holomorphic caps at the negative ends are excluded via the explicit properties of the cobordism $(C,\omega_C)$ (see \Cref{thm:strong_cobord_bourgeois} for a precise statement) and via the specific Reeb dynamics at the negative ends; this is a subtle point. 
Now, the existence of such holomorphic plane contradicts hypertightness of each connected component of the concave boundary $(C,\omega_C)$, thus concluding the proof.

\medskip

\textsc{Obstructions to fillability in arbitrary dimensions.} 
The proof of \Cref{thm:homology_injection_boundary} is mostly based on holomorphic curve techniques.

We first attach the symplectic cobordism of \cite[Section 6.1]{MNW} to the original contact manifold to obtain a symplectic fibration over $\partial V$, with $\SSS^2$ as fiber; see Lemma \ref{lemma:capping_cobord_high_dim}. This cobordism can be seen as the attaching of a symplectic handle, having two distinguished symplectic co-cores $C_\pm$. 

In the presence of a strong symplectic filling $W$, after attaching such a symplectic cobordism and obtaining a  ``capped'' filling $W_{cap}$ with boundary a fibration over a contact manifold with symplectic sphere fibers, we obtain an induced moduli space $\mathcal{M}$ of closed holomorphic spheres which ``probes'' the original filling $W$. A key point is that each of the two co-cores $C_\pm$ of the relevant handle of the capping are $J$-invariant, and every curve in the moduli space intersects it precisely once, by positivity of intersections. The same is true after considering the Gromov compactification $\overline{\mathcal{M}}$ obtained by adding strata of nodal curves. 

Assuming that there is a cycle $\sigma$ in $\partial C_\pm$, which is non-trivial in $C_\pm$ and bounds a relative cycle $b$ in the filling, this allows one to ``pull-back'' $b$ to the moduli space, and push it to a relative cycle in the co-core bounding the original $\sigma$, thus contradicting the non-triviality of $\sigma$ in $C_\pm$. 
Such pull-backs can be achieved via the theory of pseudo-cycles in the case that $W$ is semi-positive; in general, we appeal to polyfold theory \cite{HWZ_GW}  as, for example, in \cite[Section 7.2]{MNW}.
In other words, the homology of the co-core must then survive in the filling. 

\medskip

\textsc{Topology of aspherical fillings.} The proof of \Cref{thm:hom_apsherical_filling_Bour} relies on the fact that the moduli space considered above is automatically compact (see \Cref{prop:compactnessTn}). One then can puncture the spheres by removing their intersection with (small open neighbourhoods of) the co-cores, obtaining a moduli space of cylinders. Since each cylinder retracts onto any of its boundary components, this moduli space also admits a retraction to a piece of its boundary, along which the evaluation map is a diffeomorphism. Using a pull-push argument as in the proof of \Cref{thm:homology_injection_boundary}, one obtains the desired isomorphism in homology.

In the case of $S^*\TT^n$ applying some standard algebraic topology, one can then prove that any symplectically apherical filling $W$ is in fact homotopy equivalent to $D^*\TT^n$. Then, using the $s$--cobordism theorem as in \cite[Section 8]{BGZ}, one concludes that $W$ is diffeomorphic to $D^*\TT^n$. For this final step it is essential that the fundamental group of $S^*\TT^n$ is abelian. 

\medskip

\textbf{Acknowledgments.} The authors wish to thank Chris Wendl, for his hospitality when receiving the first and second authors in Berlin, and for helpful discussions with all three authors throughout the duration of this research. 
We also thank Sebastian Hensel for suggesting a simplified proof of \Cref{corAnosov}, as well as Bahar Acu, Zhengyi Zhou, Ben Filippenko, Kai Cieliebak, Sam Lisi, Patrick Massot, Klaus Niederkr\"uger, Fran Presas, Richard Siefring, Andr\'as Stipsicz and Otto van Koert for many helpful discussions and comments on preliminary versions of this paper. We are particularly grateful to Brett Parker and Urs Fuchs for pointing out a flaw in an earlier version. We further thank the  anonymous referees for very detailed and careful suggestions, which have significantly improved the results and content of this article.

The first and second authors greatly acknowledge the hospitality of the MATRIX center, where part of this research was conducted. Similarly, the second and third authors are grateful to Monash University. The second author has been supported by the grant NKFIH KKP 126683 and the European Research Council (ERC) grant No. 772479 during the preparation of this work. The third author acknowledges the support by the National Science Foundation under Grant No. DMS-1926686, during his Membership at the Institute for Advanced Study in Princeton, NJ.

\section{\textbf{The Bourgeois construction}}
\label{sec:bourgeois_constr}

Consider a closed, oriented, connected smooth manifold $V^{2n-1}$ and an open book decomposition $(B,\theta)$, together with a defining map $\Phi\co V\rightarrow \DD^2$ having each $z\in\Int(\DD^2)$ as regular value. Here, $B\subset V$ is a closed codimension-2 submanifold, $\theta=\Phi/\norm{\Phi}: V \backslash B \rightarrow \mathbb{S}^1$ is a fiber bundle, and $\Phi$ is such that $\Phi^{-1}(0)=B$.

A $1$-form $\alpha$ on $V$ is said to be \emph{adapted} to $\Phi$ if it induces a contact structure on the regular fibers of $\Phi$ and if $d\alpha$ is symplectic on the fibers of $\theta=\Phi/\norm{\Phi}$.
In particular, if $\xi$ is a contact structure on $V$ supported by $(B,\theta)$, in the sense of \cite{Gir}, then (by definition) there is such a pair $(\alpha,\Phi)$ with $\alpha$ defining $\xi$.

\begin{theorem}[Bourgeois \cite{Bo}]
	\label{thm:bourgeois}
	Consider an open book decomposition $(B,\theta)$ of $V^{2n-1}$, represented by a map $\Phi=(\Phi_1,\Phi_2)\co V \rightarrow \RR^2$ as above, and let $\alpha$ be a $1$-form adapted to $\Phi$.
	Then, $\beta \coeq \alpha + \Phi_1 d q_1 - \Phi_2 dq_2$ is a contact form on $M\coeq V\times \TT^2$, where $(q_1,q_2)$ are coordinates on $\TT^2$.
\end{theorem}
The contact form $\beta$ on $M=V\times\TT^2$ will be called \emph{Bourgeois form} associated to $(\alpha,\Phi)$ in the following. 

\begin{remark}
\label{rmk:bourgeois_form_well_def}
The contact structure on $M=V\times\TT^2$ defined by $\beta$ is actually independent, up to isotopy, on the pair $(\alpha,\Phi)$ defining the contact open book $(B,\theta)$ on $V$.
This can easily seen in the case where $n=2$, i.e.\ $V$ is $3-$dimensional, using the contractibility of the space of symplectic forms on the $2-$dimensional page. 
In the higher dimensional setting, this is discussed in detail in \cite[Section 2 and Appendix]{LMN}, using the reinterpretation of contact open books in terms of ideal Liouville structures in \cite{GirouxIdealLiouvDom}.
\end{remark}

\begin{remark}\label{rem:cover_finite}
	The contact structure determined by a Bourgeois contact form is stable up to contactomorphism under finite covers of the torus factor. 
	Indeed, up to precomposing by an automorphism of $\TT^2$, any such cover is of the form 
	$$(q_1,q_2) \longmapsto (kq_1,q_2)\text{ .}$$
	Pulling back gives a contact form $\beta_k = \alpha + k\Phi_1 d q_1 - \Phi_2 dq_2$ and a straightforward calculation shows that linear interpolation gives a family of contact forms.
\end{remark}

\medskip
\textbf{Abstract open books and Bourgeois contact structures.} 
For the proof of \Cref{thm:tight_bourgeois}, it is also useful to interpret the Bourgeois construction in abstract terms. We briefly recall here the construction in order to fix some notation.
The reader can consult for instance \cite[Section 7.3]{Gei} for further details.

Consider a Liouville domain $(\Sigma^{2n-2},\lambda)$, together with an exact symplectomorphism $\psi$ of $(\Sigma,d\lambda)$ (i.e.\ $\psi^* \lambda=\lambda - dh$, for some smooth $h:\Sigma \rightarrow \mathbb{R}^+$), fixing pointwise a neighborhood of the boundary $B\coeq \partial\Sigma$.
One can then consider the mapping torus $\Sigma_\psi$ of $(\Sigma,\psi)$, and the abstract open book
\begin{equation}
\label{eq:abstract_open_book}
V_{\Sigma,\psi}\coeq \bigslant{\left(B\times\DD^2 \sqcup \Sigma_\psi\right)}{\sim} 
\end{equation}
where $\sim$ identifies $(p,\theta)\in \partial(B\times\DD^2)$ with $[p,\theta]\in\partial \Sigma_\psi$.
One can also construct a fiberwise Liouville form $\lambda_{\psi}$ on the mapping torus $\pi_0\co\Sigma_\psi\rightarrow\SSS^1$ of $(\Sigma,\psi)$.
For large $K \gg 0$ the form $\alpha_K = K\pi_0^*d\theta + \lambda_\psi$ is contact on $\Sigma_\psi$. 
Moreover, it can be extended to a contact form on all of $V_{\Sigma,\psi}$ by
$h_1(r)\lambda_{B}+h_2(r)d\theta$ on  $B\times\DD^2$, for a well chosen pair of functions $(h_1,h_2)$, and $\lambda_B=\lambda\vert_{B}$.

We denote the resulting contact form on $V_{\Sigma,\psi}$ by $\alpha_{\Sigma,\psi}$.
The contact manifold $(V_{\Sigma,\psi},\ker(\alpha_{\Sigma,\psi}))$ will also be called an \emph{abstract contact open book}, and denoted simply with $\OBD(\Sigma,\psi)$.
Sometimes, we will also use the contact form $\alpha_{\Sigma,\psi}$.

We point out that there is a well defined map $\Phi_{\Sigma,\psi}\co V_{\Sigma,\psi}\rightarrow\DD^2$ given by extending the projection to the circle on $\Sigma_\psi$ by setting $$\Phi_{\Sigma,\psi}\vert_{B\times\DD^2}(p,r,\theta)= \rho(r)e^{i\theta} \in \DD^2,$$ for some non-decreasing function $\rho$ satisfying $\rho(r)=r$ near $0$ and $\rho(r)=1$ near $r=1$. 
Notice also that $\alpha_{\Sigma,\psi}$ is naturally adapted to $\Phi_{\Sigma,\psi}$ (as defined above).

We then denote by $\beta_{\Sigma,\psi}$ the Bourgeois form on $M_{\Sigma,\psi}\coeq V_{\Sigma,\psi}\times\TT^2$ associated to $(\alpha_{\Sigma,\psi},\Phi)$ as in \Cref{thm:bourgeois}, and by $\xi_{\Sigma,\psi}$ the contact structure it defines. Finally we let $BO(\Sigma, \psi):=(M_{\Sigma,\psi}, \xi_{\Sigma,\psi})$.

\medskip
\textbf{Hypertightness for Bourgeois Contact Forms.}
In the following sections, we will need a hypertightness criterion for $\alpha_{\Sigma,\psi}$. We first give a definition:
\begin{def}
	\label{def:poss_contr_reeb_dyn}
	Let $(V\times\TT^2,\xi)$ be a contact manifold, and $\calB$ any subset of the set of closed Reeb orbits of a contact form $\beta$. 
	We say that $\beta$ has \emph{$\TT^2$-trivial Reeb dynamics concentrated in $\calB$} if the image of every closed Reeb orbit not in $\calB$ under the projection $V\times\TT^2\rightarrow \TT^2$ is homotopically non-trivial.
\end{def}

A straightforward computation gives:
\begin{obs}\emph{\cite[Section 10.2]{Bothesis} (cf.\ \cite[Corollary 6.3]{Gi}).}
	\label{hypertightobs}
	The Bourgeois contact form $\beta_{\Sigma,\psi}$ for $\xi_{\Sigma,\psi}$ has $\TT^2$-trivial Reeb dynamics concentrated in the set $\calB$ consisting of the submanifolds $\gamma_B\times\{q\}\subset V\times\TT^2$, for all $q\in\TT^2$ and all $\gamma_B$ closed Reeb orbit of $(B,\alpha_{\Sigma,\psi}\vert_{B})$.
	If the binding $(B,\alpha_{\Sigma,\psi}\vert_{B})$ of the natural open book of $V_{\Sigma,\psi}$ admits no Reeb orbits that are contractible in $V_{\Sigma,\psi}$, then the Bourgeois contact structure $\xi_{\Sigma,\psi}$ is hypertight. 
\end{obs}
Notice that, in the $3$-dimensional case, \Cref{hypertightobs} implies that, if the binding consists of a collection of loops each having infinite order in $\pi_1(V)$, then the associated Bourgeois contact structure is hypertight.

We point out that we will not make use of \Cref{hypertightobs} in the proof of \Cref{thm:tight_bourgeois}, as we will apply it directly on a another contact form, which still defines the Bourgeois contact structure up to isotopy (see \Cref{lemma:existence_plane} below).

\medskip
\textbf{A supporting spinal open book decomposition.} We now present a geometric way of understanding the Bourgeois construction, via SOBDs (see also \Cref{sec:general_obstr_fillab} for an alternative SOBD). The notion of an SOBD, introduced in \cite{LvHMW} in dimension $3$ (see also \cite[Appendix B]{MorPhD} for a version in arbitrary dimensions), comes from the observation that an open book decomposition can be thought of as two fibrations glued together: the neighbourhood of the binding is a contact fibration over the $2$-disc, while the mapping torus piece is a Liouville fibration over the circle. More generally, an SOBD consists of a contact fibration over a general Liouville domain, glued to a Liouville fibration over a contact manifold. The first fibration is called the \emph{spine}; the second one, the \emph{paper}. The base of the spine (the \emph{vertebrae}) has boundary the base of the paper; the fibers of the paper are called the \emph{pages}. The two fibrations are glued together along an \emph{interface region}, which for our purposes we will think as a trivial product of a piece of the vertebrae (a collar neighbourhood of its boundary) and the contact fiber of the spine, and hence can be given the structure of both types of fibrations.

Let us see how this works in the case of Bourgeois manifolds. Consider $\Phi=\Phi_{\Sigma,\psi}=\rho e^{i\theta}=(\rho\cos(\theta), \rho\sin(\theta))$, a defining map for $V=V_{\Sigma,\psi}=OBD(\Sigma,\psi)$, together with the Giroux form $\alpha=\alpha_{\Sigma,\psi}$ and the associated Bourgeois form $\beta=\beta_{\Sigma,\psi}$. 
Let $\theta=\Phi/|\Phi| \co V \setminus B \to \SSS^1$ be the open book coordinate. 
From \Cref{eq:abstract_open_book}, we obtain a decomposition
$$
M=V \times \mathbb{T}^2= B \times D^*\mathbb{T}^2\cup \Sigma_\psi \times \mathbb{T}^2 \text{ ,}
$$
where we identify $D^*\mathbb{T}^2\xrightarrow{\sim}\mathbb{D}^2\times \mathbb{T}^2$ via $(q_1,p_1,q_2,p_2)\mapsto (p_1,-p_2,q_1,q_2)$.

We denote by $M_S:= B \times D^*\mathbb{T}^2$, which we call the \emph{spine}, and $M_P:=\Sigma_\psi \times \mathbb{T}^2$, the \emph{paper}. We also have the \emph{interface region} $M_I\cong B \times [-\epsilon,\epsilon]\times \mathbb{T}^3$, corresponding to the region where $M_S$ and $M_P$ glue together. Observe that we have fibrations
$$\pi_S: M_S \rightarrow D^*\mathbb{T}^2, \quad \pi_P: M_P \rightarrow S^*\mathbb{T}^2=\mathbb{T}^3,$$
where the monodromy of $\pi_P$ coincides with $\psi$ along the cotangent $\SSS^1$-direction, and is trivial along $\mathbb{T}^2$. 
The map $\pi_S$ has contact fibers and Liouville base, whereas $\pi_P$ has contact base and Liouville fibers. The interface region is a trivial product and hence fibers over $B$ or over a collar neighbourhood of $\partial D^*\mathbb{T}^2$ inside $D^*\mathbb{T}^2$. This is then an SOBD for $M$. Observe that the fibers of $\pi_P$, the \emph{pages} of the SOBD, coincide with the pages of the OBD for $V$. 
One may also view the SOBD as a fibration $\widehat{\pi}_P: M\backslash \mathcal{B}\rightarrow \mathbb{S}^1\times \mathbb{T}^2$, where we define the \emph{binding} of the SOBD as $\mathcal{B}=B \times \{0\} \times \mathbb{T}^2 \subset B \times \DD^2 \times \TT^2 = M_S$. 
This fibration has fibers which symplectically are copies of the Liouville completion of the page $\Sigma$, and has monodromy $\psi$ along the first factor, and trivial along the second one. 

The Bourgeois contact structure $\xi=\xi_{\Sigma,\psi}$ is ``supported'' by the SOBD described above, in a sense which we now describe. Via the identification $D^*\mathbb{T}^2\rightarrow\mathbb{D}^2\times \mathbb{T}^2$ above, up to isotopy of contact forms, we have
$$\beta\vert_{M_S\backslash M_I}= \lambda_B + \lambda_{std},$$
where $\lambda_{std}=p_1 dq_1+p_2dq_2$ is the standard Liouville form on $D^*\mathbb{T}^2$. In other words, $\beta\vert_{M_S\backslash M_I}$ is a split contact form, having a Liouville and a contact summand. 
Note also that on $M_S\backslash M_I$ the Reeb vector field $R_\beta$ of $\beta$ agrees with $R_B$ over the binding $\mathcal{B}$, and is transverse to the pages away from it. Similarly, up to isotopy $$\beta\vert_{M_P\backslash M_I}=\lambda_\psi + \alpha_{std},$$ where $\alpha_{std}=\cos(\theta)dq_1+\sin(\theta)dq_2$ is the standard contact form on $\mathbb{T}^3$, and so splits into a Liouville summand and a contact summand. In particular, the restriction of $d\beta$ to the pages of the SOBD is a positive symplectic form, and the Reeb vector field is transverse to the pages, agreeing with that of $\alpha_{std}$ and so tangent to the $\mathbb{T}^2$ factor. In other words, the contact structure, the contact form, as well as the Reeb dynamics of $\beta$ are ``compatible'' with the underlying geometric decomposition. This interpretation also allows us to reobtain \Cref{hypertightobs}. 

\begin{remark}
The observation that Bourgeois contact manifolds are supported by the above SOBD, in the sense described above, should be attributed to Sam Lisi.
\end{remark}

\section{\textbf{A cobordism of Bourgeois manifolds}}
\label{sec:cobordisms}

We describe a strong (actually, pseudo--Liouville) cobordism between Bourgeois contact manifolds with the same page.
Its purpose is to relate the Bourgeois manifold coming from two different monodromies to the one coming from their composition, and it will be used in the proof of \Cref{thm:tight_bourgeois}.

\subsection{\textbf{From disjoint union to composition of monodromies}}
\label{sec:stab_klukas}

Let $(\Sigma^{2n-2},d\lambda)$ be a Liouville manifold, and let $\phi$ be an exact symplectomorphism relative to the boundary.
Notice that the boundary $(B,\lambda_B)\coeq(\partial \Sigma,\lambda\vert_{\partial\Sigma})$ can naturally be seen as the ``binding" submanifold of the associated open book. 
For each $q\in\TT^2$, we also let $B_{q}$ be $B\times\{q\}\subset V_{\Sigma,\phi}\times\TT^2= M_{\Sigma,\phi}$.

The aim of this section is to give a proof of the following result (recall Figure \ref{StabKlfig}):

\begin{theorem}
	\label{thm:strong_cobord_bourgeois} 
	There is a smooth cobordism $C$ from $M_{\Sigma,\psi} \bigsqcup M_{\Sigma,\phi}$ to $M_{\Sigma,\psi \circ \phi}$. 
	This cobordism is smoothly a product $C_0\times\TT^2$, where $C_0$ is a smooth cobordism  from $V_{\Sigma,\psi} \bigsqcup V_{\Sigma,\phi}$ to $V_{\Sigma,\psi \circ \phi}$. 
	Moreover, there is a symplectic form $\omega_C$ on $C$ which satisfies the following properties:
	\begin{enumerate}
	    \item \label{item:sympl_form_1} $\omega_C$ admits local Liouville forms $\lambda_+$ and $\lambda_-$ near $M_{\Sigma,\psi \circ \phi}$ and $M_{\Sigma,\psi} \bigsqcup M_{\Sigma,\phi}$ respectively, satisfying:
	    \begin{enumerate}
	        \item \label{item:local_liouv_1} $\lambda_+$ restricts on $M_{\Sigma,\psi \circ \phi}$ to the Bourgeois contact form $\beta_{\Sigma,\psi \circ \phi}$,
	        
	        \item \label{item:local_liouv_2} $\lambda_-$ restricts on $M_{\Sigma,\psi}\sqcup M_{\Sigma,\phi}$ to the Bourgeois contact forms $\beta_{\Sigma,\psi}$ and $\beta_{\Sigma,\phi}$ respectively;
	        in particular, $\lambda_-$ has (on each connected component) $\TT^2$-trivial Reeb dynamics concentrated in $\{B_q\}_{q\in\TT^2}$; 
	
	    \end{enumerate}
	    \item \label{item:sympl_form_2} $\omega_C$ admits a global primitive $\nu$ which coincides with $\lambda_+$ at the convex boundary and such that $\nu\vert_{B_q}=\lambda_-\vert_{B_q}$ for each $B_q\subset M_{\Sigma,\psi} \bigsqcup M_{\Sigma,\phi}$.
	\end{enumerate}
\end{theorem}

\noindent \Cref{item:sympl_form_1} means in particular that $(C,\omega_C)$ is a strong symplectic cobordism with convex boundary $\BO(\Sigma,\phi\circ\psi)$ and concave boundary $\BO(\Sigma,\phi)\sqcup \BO(\Sigma,\psi)$.
Notice however that we do {\em not} claim that the global $1$-form $\nu$ in \Cref{item:sympl_form_2} defines a contact structure at the concave boundary;
in other words, the cobordism we give is not claimed to be Liouville, but just pseudo-Liouville (as defined in the introduction).
Lastly, we point out that \Cref{thm:strong_cobord_bourgeois} can be thought of as a ``stabilized'' version of \cite[Proposition 8.3]{Advek} and \cite[Theorem 1]{Klu12};
in fact, smoothly (but not symplectically), the cobordism $C$ is just the product of the cobordism from \cite{Advek,Klu12} with $\mathbb{T}^2$.
For the reader's convenience, we start by giving a topological description of the cobordism $C$ as obtained by gluing two ``cobordisms with corners'', $\Cbot$ and $\Ctop$, and then describe the symplectic structures on these pieces in more detail.

\medskip 

\textbf{A topological description of $C$.}
Let $P$ be the pair of pants, i.e.\ the surface of genus $0$ and with $3$ boundary components. 
We view $P$ as embedded in $\RR^2$ as the (closed) disc $D_+$ of radius $1$ with two smaller (open) disjoint discs $D_{-,1}$ and $D_{-,2}$, both of radius $\epsilon$ and centered at $-1/2$ and $+1/2$ respectively, removed from it. 
In cobordisms terms, $P$ is seen as a smooth cobordism with concave boundary $\partial D_{-,1} \sqcup \partial D_{-,2}$ and convex boundary $\partial D_{+}$.

Consider then the fiber bundle $\pi_E\colon E\to P$, with fiber the page $\Sigma$, over the pair of pants $P$, where the monodromies along the two negative boundary components are given by $\phi$ and $\psi$ respectively, and by their composition along the positive one.
This can be realized for instance as follows:
consider $\gamma_1$ and $\gamma_2$ disjoint arcs in $P$ joining respectively $\partial D_{-,1}$ and $\partial D_{-,2}$ to $\partial D_+$, then cut $\Sigma\times P$ along $\Sigma\times \gamma_1$ and $\Sigma\times \gamma_2$ and glue them back respectively via $\phi\times\Id_{\gamma_1}$ and $\psi\times\Id_{\gamma_2}$. 

\begin{figure}[t]
	  \centering
	  \includegraphics[width=0.4\linewidth]{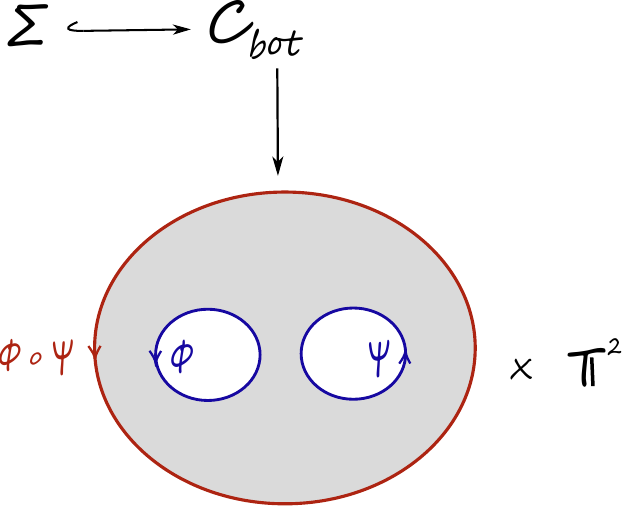}
	  \caption{Topological picture of $\Cbot$.}
	  \label{FigBottomPartCob}
\end{figure}

The desired ``bottom piece'' $\Cbot$ of $C$ is then the cobordism with corners $E\times \TT^2$, which inherits a fibration $\pi=\pi_E\times id_{\TT^2}: C_{bot}\rightarrow P\times \mathbb{T}^2$; cf.\ \Cref{FigBottomPartCob}.
Notice that $\Cbot$ has the following distinguished boundary pieces:
\begin{itemize}
    \item[$\bullet$] $\partial_{\pm}\Cbot$ given by the preimage by $\pi$ of $\partial D_+\times \mathbb{T}^2$ and $\partial D_{-,1}\sqcup \partial D_{-,2}\times \mathbb{T}^2$, respectively;
    \item[$\bullet$] $\partial_0 \Cbot$ given by $\partial \Cbot \setminus \pi^{-1}(\partial P\times \mathbb{T}^2)$.
\end{itemize}
Furthermore, since $\phi$ and $\psi$ are both identity near $\partial \Sigma$, one simply has that $\partial_0 \Cbot = \partial \Sigma \times P \times \mathbb T^2$.

\begin{figure}[t]
	  \centering
	  \includegraphics[width=0.5\linewidth]{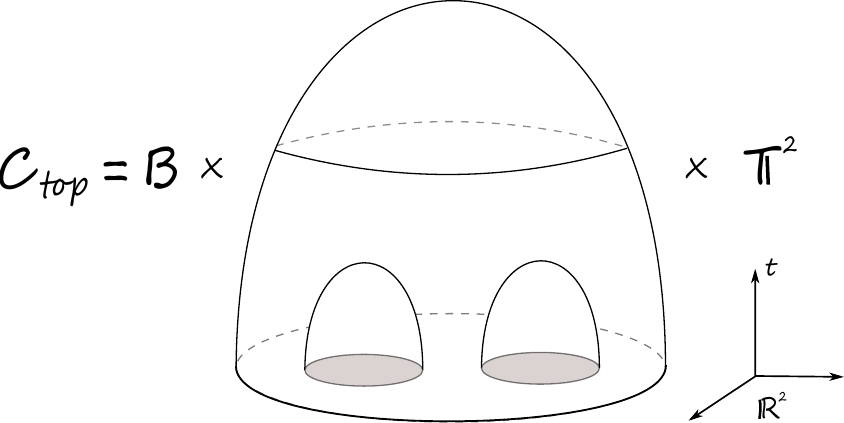}
	  \caption{\normalsize Topological picture of the cobordism (with corners) $\Ctop\subset B\times [0,+\infty)\times \RR^2\times\TT^2$.
	  In the proof of \Cref{thm:strong_cobord_bourgeois}, $\RR^2\times\TT^2$ becomes $T^*\TT^2$ and  $[0,+\infty)\times B$ becomes the symplectization of $B$.}
	  \label{FigTopPartCob}
\end{figure}

The ``top piece'' $\Ctop$ of $C$ is given as follows.
Consider $P\subset \RR^2$ as subset of $\{0\}\times\RR^2\subset\RR^3$. 
Let also $S_+$, $S_{-,1}$ and $S_{-,2}$ be three hemispheres in $[0,+\infty)\times\RR^2$ centered respectively at $0$, $-1/2$ and $+1/2$, and of radius $1$, $\epsilon$ and $\epsilon$ respectively.
Then, $P$ together with these three hemispheres bounds a compact region $R\subset \RR^3$ with piece-wise smooth boundary.
The desired $\Ctop$ is then just given by $\Ctop= B \times R \times \TT^2$, where $B=\partial \Sigma$; cf.\ \Cref{FigTopPartCob}.
As was the case with $\Cbot$, the manifold $\Ctop$ also has the distinguished boundary pieces:
\begin{itemize}
    \item[$\bullet$] $\partial_{+}\Ctop$ given by $ B \times S_{+} \times \TT^2$;
    \item[$\bullet$] $\partial_{-}\Ctop$ given by $B\times(S_{-,1}\sqcup S_{-,2})  \times \TT^2$;
    \item[$\bullet$] $\partial_0 \Ctop$ given by $ B \times P \times \TT^2$.
\end{itemize}

\noindent The desired cobordism $C$ is then topologically just obtained by gluing $\Ctop$ with $\Cbot$ along $\partial_0 \Ctop = B \times P \times \TT^2$ and $\partial_0 \Cbot= B \times P \times \TT^2$ (via the natural identification).

\medskip
\medskip 

We now proceed to discuss the symplectic side of the construction, i.e.\ to give a detailed proof of \Cref{thm:strong_cobord_bourgeois}.
We will follow closely, with some adaptations, the proof given in \cite{Klu12}.

\medskip

\textbf{A toroidal pair of pants cobordism in dimension $4$.} 
In order to make the topological sketch above into a detailed proof taking care of the symplectic data, we need to utilise a symplectic counterpart of the product of the pair of pants $P$ and $\TT^2$ which has the right structure at the boundary. 
In other words, we want a strong $4$-dimensional symplectic cobordism with concave end $(S^*\TT^2,\xi_{std})\sqcup (S^*\TT^2,\xi_{std})$ and convex end $(S^*\TT^2,\xi_{std})$, and having the topology of $P\times \TT^2$. 
To this end, consider the unit disc cotangent bundle $D^*\TT^2$ of $\TT^2$, together with its standard symplectic structure $\omegastd= d\lambdastd$, where in coordinates $\lambdastd=p_1dq_1+p_2dq_2$ is the standard Liouville form. 
To be precise, we need to work with scalar multiples $K\omegastd$ and $K\lambdastd$, where $K$ is a positive real constant that will be determined later on in the proof. We also denote by $X$ the Liouville vector field $p_1\partial_{p_1}+ p_2\partial_{p_2}$.

Consider the submanifold $D_\epsilon^*\TT^2$ of $D^*\TT^2$ made of those covectors of norm less than a certain $\epsilon<1/10$, and denote by $j_{\pm}\co D_\epsilon^*\TT^2\rightarrow D^*\TT^2$ the symplectomorphisms 
\begin{equation*}
    (p_1,q_1,p_2,q_2)\stackrel{j_\pm}{\longmapsto}(p_1\pm 1/2,q_1,p_2,q_2) \text{.}
\end{equation*}
For ease of notation, we consider the inclusion $$j=j_-\sqcup j_+\colon (D_{\epsilon}^*\TT^2, K\omegastd)\sqcup (D_{\epsilon}^*\TT^2, K\omegastd)\rightarrow (D^*\TT^2, K\omegastd).$$
Then, the desired cobordism is $(Q,\omega_Q)\coloneqq (D^*\TT^2\setminus j(D_\epsilon^*\TT^2), K\omegastd)$. Topologically this is just a product of the torus with a pair of pants $P \subset \mathbb{R}^2$ in the $(p_1,p_2)$-plane.

The Liouville field on a neighbourhood of the convex boundary is just given by $X$, whereas the one near the concave boundary is given by $j_*X$, which is the vector field $(p_1\mp 1/2)\partial_{p_1}+p_2\partial_{p_2}$ on the image of $j_\pm$ respectively.

Lastly, we consider an auxiliary smooth function $f\colon T^*\TT^2 \rightarrow \RR$, which depends only on $p_1,p_2$, and satisfies:
\begin{enumerate}
	\item \label{item:cond_f_1} $f=p_1^2+p_2^2$ on $T^*\TT^2\setminus D^*\TT^2$,
	\item \label{item:cond_f_2} $f=(p_1\mp1/2)^2+p_2^2$ on the image of $j_\pm$ respectively,
	\item $\epsilon^2 < f < 1$ on the interior of $Q$.
\end{enumerate} 
Notice that \Cref{item:cond_f_1} implies in particular that $df(X)>0$ in a neighborhood of $\partial (D^*\TT^2)$. 
Similarly, \Cref{item:cond_f_2} implies that $df(j_*X)>0$ on the image of $j$, except on $j(\{0\}\times\TT^2\sqcup \{0\}\times\TT^2)$, where it vanishes.

\medskip

\textbf{Description of the symplectic cobordism.}
Consider the bottom piece $\Cbot$ as in the topological sketch above. 
We now view $\Cbot$ as a fiber bundle with fibers $\Sigma$ over $Q=P\times \TT^2 \subset T^*\TT^2$, and want to prove that it admits a fiberwise Liouville form $\lambda_0$ which agrees with $e^t\lambda_B$ near $B\times Q\subset \Cbot$.
Here, $e^t\lambda_B$ is the normal form of $\lambda$ on $\Sigma$ on a sufficiently small neighborhood $(-\delta,0]\times B$, with $t\in(-\delta,0]$, of its contact boundary $B$ given by the (globally defined) Liouville vector field $Y$ for $\lambda$, which is outwards pointing along the boundary $\partial \Sigma = B$ of the page $\Sigma$.

Recall that $\Cbot$ is obtained from $\Sigma \times T^*\TT^2$ by a cut and paste procedure along $\Sigma\times(\gamma_1\cup\gamma_2)\times\TT^2$, as described in the topological sketch above. 
Recall also that $\phi$ and $\psi$ are exact symplectomorphisms of the page $(\Sigma,\lambda)$, equal to the identity on a neighborhood of $B$; let then $h_1$ and $h_2$ be respectively the functions on $\Sigma$, constant near the boundary, such that $\phi^*\lambda = \lambda - dh_1$ and  $\psi^*\lambda = \lambda - dh_2$.
Taking normal coordinates $r_i\in(-\varepsilon,\varepsilon)$ to $\gamma_i$ inside $P$, one can then consider smooth cutoff functions $\rho_i\colon (-\varepsilon,\varepsilon)\to [0,1]$ equal to $0$ near $r_i=-\varepsilon$ and to $1$ near $r_i=0$.
It then follows that, in the cut and paste procedure to obtain $\Cbot$ from $\Sigma\times T^*\TT^2$, the form $\lambda-d(\rho_i h_i)$ on $\Sigma \times(P\setminus (\gamma_1 \cup\gamma_2)) \times \TT^2$ glues well under the identifications $\phi\times\Id_{\gamma_1}$ and $\psi\times\Id_{\gamma_2}$ along respectively $\Sigma \times \gamma_1\times\TT^2$ and $\Sigma \times \gamma_2\times\TT^2$.
The result is hence a fiberwise Liouville form $\lambda_0$ on the fiber bundle $\pi\colon\Cbot\to Q$, which is equal to $e^t\lambda_B$ near $B\times Q$ as the $h_i$'s are constant there.

We also point out that, for very large $K \gg 0$, the form $\Omega = d\lambda_{0}+ K\pi^* d\lambda_{std} $ on the total space $\Cbot$ of $\pi\colon \Cbot\to Q$ is an exact symplectic form with primitive $\nu = \lambda_{0}+ K\pi^* \lambda_{std}$. 
Here, abusing notation sightly, we let $\lambda_{std}$ denote  the restriction of the canonical $1$-form under the inclusion $Q\subset T^*\TT^2$. This then gives the desired symplectic structure on $\Cbot$; it only remains to describe the symplectic form on $\Ctop$.

For this, let $\tau$ be a function on $[0,\infty)$ vanishing at $0$ with all its derivatives, and strictly monotone increasing on $(0,\infty)$.
We then identify the $\Ctop$ of the topological sketch above with the following set: 
%We then consider the manifolds with corners $\Cbot\coeq E$ and
\begin{equation*}
%\resizebox{.99\hsize}{!}{
\Ctop = \left\{(t,b,p,q)\in [0,\infty)\times B\times T^*\TT^2 \,\left\vert\, \tau(t)^2+f(p,q)^2 \geq \epsilon^2, \tau(t)^2 + \norm{p}^2 \leq 1 \right. \right\} \text{ .}
%}
\end{equation*}
The advantage of this identification is that $\Ctop$ inherits a natural Liouville structure $e^t\lambda_B+K\lambda_{std}$ from the symplectization $([0,\infty)\times B , d(e^t\lambda_B))$ and the cotangent structure $(T^*\TT^2,\lambda_{std})$.

Now, $\Cbot$ and $\Ctop$ can be naturally glued along the subsets $B\times Q\subset \partial \Cbot$ and $\{\tau=0\}=\{0\}\times B\times Q\subset \Ctop$ (corresponding to $\partial_0\Ctop$ and $\partial_0\Cbot$ in the topological sketch).
What's more, this gluing is compatible with the symplectic form $\Omega=d\lambda_{0}+ K\pi^* d\lambda_{std} $  on $\Cbot$ and $d(e^t\lambda_B)+Kd\lambda_{std}$ on $\Ctop$.  
We thus obtain a symplectic cobordism $(C,\omega_C)$, and we now proceed to check that it satisfies the required properties. 
We will only do this for the negative ends, as the case of the positive end is entirely analogous.
%\begin{figure}[t]
%	  \centering
%	  \includegraphics[width=0.35\linewidth]{bottom_cobord_2.pdf}
%	  \caption{Picture of $\Cbot$.}
%	  \label{FigBottomPartCob}
%\end{figure}

\paragraph{Proof of \Cref{thm:strong_cobord_bourgeois}}

Notice first that by construction $C = C_0 \times \TT^2$ is topologically a product. 
Furthermore, each boundary component of $C$ determine $\TT^2$-stabilised open books, which we think of as supporting SOBD's for the Bourgeois contact structures. Moreover, these SOBD's coincide, at least topologically, with those used to define the Bourgeois contact structure.
From a symplectic point of view, we point out that the fiberwise Liouville form $\lambda_0$ on $\Cbot\to Q$ glues well to the $1$-form $e^t\lambda_B$ on $\Ctop$ inherited from the ambient space $[0,\infty)\times B \times T^*\TT^2$ to give a global $1$-form denoted $\widehat{\lambda}_0$, which is a pullback of a form on $C_0$. 

We now describe the primitives near the boundary in the bottom part $\Cbot$.
Consider the primitive $j_*\lambda_{std}$ of the standard symplectic form near the concave boundary of $Q$, which in coordinates is just $p_1dq_1+p_2dq_2 \mp 1/2 dq_1 $. 
The Liouville vector field $Z$ associated to the local primitive $\lambda_- = \lambda_{0} + Kj_*\lambda_{std}$ is then defined on a neighborhood of the negative boundary components $\partial_- \Cbot$ of $\Cbot$.
More precisely, an explicit computation shows that the vector field $Z$ is, in the explicit pre-glued version $\Sigma \times (P\setminus (\gamma_1\cup\gamma_2))\times\TT^2$ of $\Cbot$, just $Y + Y' + j_*X + X'$, where $Y$ is the Liouville form associated to $\lambda$ on $\Sigma$, $Y'$ is the $d\lambda-$dual of $\rho dh$ (hence tangent to the $\Sigma$ factor), and $X'$ is the $d\lambda_{std}-$dual of $\frac{1}{K} h d\rho$; here we define $\rho$ to be $\rho_i$, respectively on each of the two negative ends. 
Hence, the vector field $X'$ is  tangent to the second factor of $\Sigma \times (P\setminus (\gamma_1\cup\gamma_2))\times\TT^2$, and has no $\partial_p$ component, using cotangent coordinates $(p,q)\in T^*\TT^2$.
In particular, the projection of $Z$ on $\Cbot$ to $Q$ via the fiber bundle map has the same $\partial_p$ component as $j_*X$ near $\partial_- Q$, and is hence transverse to it and inwards pointing.

Now, via the natural (orientation preserving) diffeomorphisms 
\begin{equation}
\label{eqn:diffeo_stab_klukas_proof}
\begin{split}
\SSS^1_\epsilon\times\TT^2 & \to S_\epsilon^*\TT^2 \subset D^*_\epsilon\TT^2 \\
(\varphi,x,y) & \to (p_1=-\cos(\varphi) \pm 1/2,q_1=x,p_2=\sin(\varphi),q_2=y),
\end{split}
\end{equation}
where $\SSS^1_\epsilon$ is the circle in $\RR^2$ of radius $\epsilon$,
the restriction of the associated Liouville form on $\Cbot$ to the boundary is naturally pulled back to the Bourgeois contact form on the paper region of the SOBD, for a suitable choice of adapted contact form $\alpha$ and open book map $\Phi$ (cf.\ \Cref{rmk:bourgeois_form_well_def}). %, and so is adapted to the Bourgeois SOBD in that region (in the sense of \Cref{sec:bourgeois_constr}).

In the top part $\Ctop$ of the cobordism the Liouville vector field is just the linear combination $\partial_t + r\partial_r$, where $r$ is the radial parameter $|p \mp 1/2|$ at the respective boundary components. 
This is then also transverse to the boundary. 

By parametrizing the interior hemispherical caps using the canonical coordinates coming from the cotangent bundle of the torus, we obtain the following coordinate description on $\Ctop$:
\begin{equation*}
    	\lambda_{-}= e^{\sqrt{\epsilon^2-\norm{r}^2}}\lambda_B + K (p_1\mp 1/2) dq_1 + K p_2 dq_2 \text{ .}
\end{equation*}

Now, via the (orientation preserving) diffeomorphisms $$\DD_\epsilon^2\times\TT^2\to D_\epsilon^*\TT^2$$ $$(r,\varphi,x,y)\to (p_1=-r\cos(\varphi) \pm 1/2,q_1=x,p_2=r\sin(\varphi),q_2=y),$$
naturally extending \Cref{eqn:diffeo_stab_klukas_proof}, where $\DD_\epsilon^2$ is here the disc in $\RR^2$ centered at the origin of radius $\epsilon$, $\lambda_-$ then pulls back on each negative boundary component of $\Ctop$ to the Bourgeois contact form as in \Cref{thm:bourgeois}, for a suitable choice of adapted contact form $\alpha$ and open book map $\Phi$, which can be taken compatible with that coming from the bottom boundary part. 
This concludes the proof of \Cref{item:local_liouv_2}.

Finally, notice that there is a global Liouville primitive $\nu$ for $\omega_C$ defined on all $C$, given by gluing $\lambda_{0} + K\pi^*\lambdastd $ on $\Cbot=E\times \TT^2$ and $e^t\lambda_B + K \lambdastd$ on $\Ctop \subset [0,+\infty)\times B\times T^*\TT^2$. 
By looking at the explicit expression for $\lambda_-$ along $\Ctop$, we see that $\nu$ coincides with $\lambda_-$ on the subsets of the form $B_q$ defined before the statement of \Cref{thm:strong_cobord_bourgeois}. 
Indeed, these are given by $B_q=B\times \{x_\pm\}\times \{q\}$, where $x_\pm=(t=\epsilon,p_1=\pm 1/2,p_2=0)\in [0,+\infty)\times \mathbb{R}^2$ is the ``origin'' in the hemispheres $S_{-,1}$, $S_{-,2}$, i.e.\ the points of $S_{-,1}$ and $S_{-,2}$ with maximal $t$ value. 
There, we just have $\lambda_-\vert_{B_q}=e^\epsilon \lambda_B=\nu\vert_{B_q}$. This concludes the proof. $\hfill \square$
\section{\textbf{Factorizing the monodromy}}\label{sec:fact_mon}

Let $\Sigma$ denote a connected orientable surface with boundary. We will denote the mapping class group as $MCG(\Sigma)$, which is defined to be the set of isotopy classes of orientation preserving diffeomorphisms of $\Sigma$; note that we do not require these diffeomorphisms to fix the boundary components. This group is naturally isomorphic to the group of isotopy classes of diffeomorphisms of the corresponding punctured surface. One may also consider $MCG(\Sigma,\partial \Sigma)$ of mapping classes fixing the boundary, and there is a natural forgetful map
$MCG(\Sigma,\partial \Sigma) \rightarrow MCG(\Sigma)$
whose kernel is generated by boundary parallel Dehn twists.

We will refer to a surface as \emph{sporadic} if it is either a disc, an annulus or a pair of pants. These cases correspond to the mapping class group being virtually abelian. The aim of this section is to prove the following:
\begin{lemma}[Factorization Lemma]
	\label{Factlemma}
	Let $\phi\in MCG(\Sigma,\partial \Sigma)$ for a non-sporadic surface $\Sigma$. Then $\phi$ can be factored as $\phi=\phi_1 \circ \phi_2$, where, for each $i=1,2$, $\phi_i$ is such that each connected component of the binding of $V_i\coeq OBD(\Sigma,\phi_i)$ has infinite order in $\pi_1(V_i)$. 
\end{lemma}

A direct consequence of \Cref{Factlemma} and \Cref{hypertightobs} is the following:
\begin{corollary}\label{hypertightfactorization}
	Let $\phi$ be a mapping class of a compact, orientable, non-sporadic surface $\Sigma$ with boundary. Then $\phi$ can be factored as $\phi=\phi_1 \circ \phi_2$, with $\phi_1,\phi_2$ such that the Bourgeois contact manifolds $BO(\Sigma,\phi_1)$ and $BO(\Sigma,\phi_2)$ are hypertight.
\end{corollary}

In order to prove \Cref{Factlemma}, we start by recalling some results from geometric group theory and 3-dimensional hyperbolic geometry in \Cref{sec:geom_group_theory,sec:hyp_geom} respectively. 
The proof is then given in \Cref{sec:fact_lemma}.

\subsection{\textbf{Dynamics of pseudo-Anosovs acting on the curve graph}}\label{sec:geom_group_theory}
Let $\Sigma$ denote a compact, connected orientable surface (possibly with boundary). 
We recall that, by the Nielsen-Thurston classification theorem (see for instance \cite[Theorem 13.2]{FM}), every element in $MCG(\Sigma)$ or $MCG(\Sigma,\partial \Sigma)$ is either pseudo-Anosov, reducible or of finite order. 
This characterization can also be extracted from the associated action on the \emph{curve graph}, denoted $\mathcal{C}(\Sigma)$. 
This is the graph whose vertices are isotopy classes of essential, non boundary-parallel, simple closed curves on $\Sigma$ so that there is an edge between two vertices if the corresponding curves can be made disjoint via isotopy. 
Endowing the edges to have length $1$ we thus obtain a metric space on which the mapping class group acts via isometries. 

The curve graph was famously shown to be 
$\delta$-hyperbolic in the sense of Gromov by Masur-Minsky \cite{Masur_Minsky} in the case that $\Sigma$ is neither the torus nor sporadic, which has many important implications. In particular, one has the following alternative description using this action that is again due to Masur-Minsky \cite[Theorem 4.6]{Masur_Minsky}. 
\begin{theorem}[Masur-Minsky]
\label{thm_hyp_Anosov}
	Let $\phi \ne Id$ be an arbitrary mapping class on a non-sporadic compact, orientable, surface $\Sigma$ with non-empty boundary. Then, we have the following trichotomy:
	\begin{itemize}
	    \item[$\bullet$] (Finite Order): The action of $\phi$ on $\mathcal{C}(\Sigma)$ has finite orbits, but no fixed point.
	    \item[$\bullet$] (Reducible):  $\phi$ has a fixed point in  $\mathcal{C}(\Sigma)$.
	    \item[$\bullet$] (Pseudo-Anosov):  $\phi$ has no finite orbits and in fact acts \textbf{hyperbolically} on $\mathcal{C}(\Sigma)$.
	\end{itemize}

\end{theorem}
\noindent In the following we will use some basic facts about Gromov-hyperbolic spaces and their boundaries at infinity, which can for example be found in the book of Bridson-Haefliger \cite[Chapter 9]{Bridson_Haef}. 

\begin{figure}[t]
    \centering
    \includegraphics[width=0.7 \linewidth]{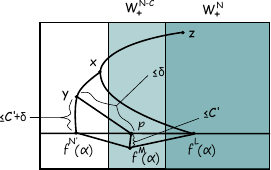}
    \caption{The triangle in the proof of \Cref{lemma:nearest_point_proj}.}
    \label{fig:pseudo_anosov_proof_1}
\end{figure}

We now collect some basic consequences of the fact that a pseudo-Anosov mapping class $f$ acts hyperbolically on the curve graph. 
First we note that, for any choice of isotopy class $\alpha$, which corresponds to a vertex in the curve graph $X =\mathcal{C}(\Sigma)$, the bi-infinite orbit $ \beta= (f^n(\alpha))_{n \in \mathbb{Z}}$ is a {\em quasi-geodesic}. Moreover, such a quasi-geodesic determines two distinct points $p_{\pm}$ on the Gromov boundary $\partial_\infty X$ of the curve graph, which correspond to the fixed points of the induced action on the boundary and form a repelling/attracting pair. Finally, the action on the boundary has {\em north-south dynamics}. More precisely, for \emph{any} neighbourhoods $U_{\pm}$ around $p_\pm \in \partial_\infty X$ there is some  (large) $N$ so that 
$$f^N(\partial_\infty X \setminus U_-) \subset U_+\textrm{ and } f^{-N}(\partial_\infty X  \setminus U_+)  \subset U_- . $$
In fact, this north-south dynamics can also be seen using half spaces associated to the quasi-geodesic $\beta = (f^n(\alpha))_{n\in \mathbb Z}$, as described presently. 
First of all the nearest point projection $\Pi_\beta: X \to \beta$, mapping an arbitrary point to a point in $\beta$ which is nearest to it, is well-defined in the coarse sense, meaning that any two nearest point projections are a uniformly bounded distance apart. We fix such a choice from now on; what follows is independent of this choice up to constants. We decompose $\beta = (f^n(\alpha))_{n\le 0} \cup (f^n(\alpha))_{n\ge 0} = \beta_+ \cup \beta_-$ into two (quasi-)rays and define half spaces
$$H_\pm = \Pi_\beta^{-1}(\beta_\pm). $$
Taking closures in $X \cup \partial_\infty X$, the north-south dynamics implies that
$$\bigcap_{n \ge 0} f^{ n}\left( \overline{H}_+\right) = p_+ \textrm{ and }\bigcap_{n \ge 0} f^{ -n} \left(\overline{H}_-\right) = p_-  \ .$$
For convenience we write $W^{N}_{\pm} = f^{\pm N}\left(H_\pm\right)$.
Note that $W^{N}_{+}$ is precisely the preimage of $\{f^n(\alpha)\}_{n \ge N}$ under the projection and similarly for $W^{N}_{-}$. In particular, $W_+^N$ and $W_-^M$ are disjoint for all $N,M \ge 1$. 
Moreover, these half spaces satisfy
the following property:
\begin{lemma}
\label{lemma:nearest_point_proj}
There is a constant $C = C(\delta,\alpha)$, depending only on the hyperbolicity constant $\delta$ and the choice of curve $\alpha$, such that every geodesic joining a point $z \in W^{N}_{\pm}$ to its nearest point projection lies in $W^{N- C}_{\pm}$.
\end{lemma}
\begin{proof}
We refer to \Cref{fig:pseudo_anosov_proof_1} for this proof.
Let $z\in W_{+}^{N}$ have nearest point projection $f^L(\alpha)$ and let $[z,f^L(\alpha)]$ be a geodesic path from $z$ to $f^L(\alpha)$. Assume that the nearest point projection $f^{N'}(\alpha)$ of some $x\in [z,f^L(\alpha)]$ sits outside $W_{+}^{N-C}$, i.e.\ $N'\leq N-C$.

As $(f^n(\alpha))_n$ is a quasi-geodesic, there is some $C'>0$ so that, for any $N-C<M<N$, there is a point $p$ on the geodesic segment $[f^{N'}(\alpha),f^L(\alpha)]$ which is at distance $C'$ from $f^{M}(\alpha)$.
Now, by the $\delta-$slim property applied to the geodesic triangle $\{x,f^{N'}(\alpha),f^L(\alpha)\}$, there is a point $y$ either on the geodesic segment $[x,f^{N'}(\alpha)]$ or on $[x, f^L(\alpha)]$ with distance at most $\delta$ from $p$.
Hence, $y$ is also at distance at most $C'+\delta$ from $f^M(\alpha)$.

Now, let us assume that $y\in [x,f^{N'}(\alpha)]$. 
Then, it follows that $[y,f^{N'}(\alpha)]$ has length at most $C'+\delta$, otherwise there would be a piecewise geodesic from $x$ to $f^M(\alpha)$ shorter than the geodesic segment $[x,f^{N'}(\alpha)]$, contradicting the fact that $f^{N'}(\alpha)$ is the nearest point projection of $x$.
But then the geodesic triangle $\{y,f^{N'}(\alpha),f^M(\alpha)\}$ violates the triangle inequality if $C$ is taken big enough (this value only depends on $\delta$ and $\alpha$, and not on the point $y$), as the two sides containing $y$ are both of length at most $C'+\delta$, while the side $[f^{N'}(\alpha),f^M(\alpha)]$ has length growing linearly in $M-N'\geq C$.

If $y\in [x,f^{L}(\alpha)]$, the analogous argument shows that the geodesic triangle $\{y,f^{L}(\alpha),f^M(\alpha)\}$ violates the triangle inequality if $C$ is big enough.
This then proves the proposition, provided $C=C(\delta,\alpha)$ is taken to be large.
\end{proof}
 
%that a geodesic joining a point $x \in W^{N}_{+}$ to its nearest point projection lies in $W^{N-C}_{\pm}$, for some constant $C = C(\delta,\alpha)$ depending only on the hyperbolicity constant $\delta$ and the choice of curve $\alpha$. 
%Another way of putting this is that any geodesic $\gamma$ with end points in  $W^{N}_{\pm}$ stays within $W^{N\mp C}_{\pm}$ for some constant as above.

We now prove the following, which only uses facts about the dynamics of isometries acting hyperbolically on  $\delta$-hyperbolic spaces. We are very thankful to Sebastian Hensel for pointing this out to us. 

\begin{proposition}\label{corAnosov}
	Let $\phi$ be a fixed mapping class on a non-sporadic compact, orientable, surface $\Sigma$ with non-empty boundary. Then, there exists a pseudo-Anosov map $f$ such that, for sufficiently large $k$, the mapping class $f^k\phi \in MCG(\Sigma,\partial \Sigma)$ is pseudo-Anosov.
\end{proposition}

\begin{proof}
Let $\phi$ be our given mapping class. 
We first claim that there is a pseudo-Anosov map $f$ whose set of fixed points $\{q_\pm\}$ at infinity is mapped to a disjoint set under $\phi$. 
To see this, first recall that, since $\Sigma$ is non-sporadic, there exist at least two pseudo-Anosov mapping classes $\gamma ,\eta \in MCG(\Sigma,\partial \Sigma)$ with disjoint sets of fixed points at infinity $\{p_\pm^\gamma\}$ and $\{p_\pm^\eta\}$\footnote{This follows, for example, from the fact that the action on the boundary of the curve complex is WPD in the sense of \cite[Proposition 6]{BF02}.}.
Now, if for instance $\phi$ mapped the stable fixed point $p_+^\gamma \in \partial_\infty X$ of $\gamma$ to its unstable fixed point $p_-^\gamma$ (the other case is similar), then we can consider the conjugate $f = \gamma^N\eta \gamma^{-N}$ for large $N$ so that both fixed points $q_\pm$ of $f$ lie in a small neighbourhood of the attracting fixed point $p_+^\gamma$ of $\gamma$, which is then mapped under $\phi$ to a small neighbourhood of $p_-^\gamma$. In particular, the set $\{\phi(q_\pm)\}$ is disjoint from $\{q_\pm\}$ as claimed.

Consider now such a pseudo-Anosov $f$, and let $W_\pm^N = f^{\pm N}(H_\pm)$ as in the notation described before the statement of \Cref{corAnosov}. Note that these subsets are defined with respect to some {\em a priori} fixed quasi-geodesic $(f^n(\alpha))_{n\in\ZZ}$. 
Observe that, in view of our assumption on $f$, we can use the north-south dynamics to assume that \emph{both} $\phi^{-1}(W_-^N)$ and $\phi(W_+^N)$ are disjoint from $W_-^N$, for all large $N \gg 0$.
%Now for a sufficiently large power, the dynamics of $f^k$ are highly \textcolor{red}{does highly mean something precise here?} contracting/repelling. 
In particular, for any given (large) $M$ we can take $k$ so that $g_k = f^k \phi$ satisfies
$$ g_k(W_+^N) \subset W_+^{N + 3M}.$$
%Furthermore, the distance between $\psi^N(\alpha)$ and $W_+^{N + M}$, which is the same as the distance of $\alpha$ to $W_+^{M}$, goes to infinity as $M\to \infty$.

 Now suppose $g_k$ is not pseudo-Anosov.
 In particular, up to taking powers, it has a fixed point, say $\alpha'$.
 We claim that such a fixed point must belong to $W_+^N$. 
Assume not. Then consider the geodesic joining $\alpha'$ to $f^{N+M}(\alpha)$. This has one end point fixed by $g_k$, while the other gets mapped to $g_k(f^{N+M}(\alpha))\in W_+^{N+3M}$ under $g_k$. Joining $g_k(f^{N+M}(\alpha))$ by a geodesic to its nearest point projection $f^{L}(\alpha)$, where $L \ge N+3M$, we obtain a geodesic quadrilateral with vertices $g_k(f^{N+M}(\alpha)),f^{L}(\alpha),$ $f^{N+M}(\alpha),\alpha'$ that has two sides of equal length $D = d(\alpha', f^{N+M}(\alpha))$ at the vertex $\alpha'$.
 See \Cref{fig:pseudo_anosov_proof_2}.
 
 \begin{figure}[t]
    \centering
    \includegraphics[width=0.7 \linewidth]{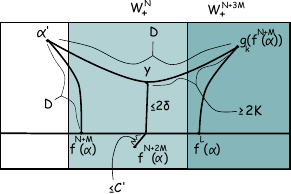}
    \caption{The quadrilateral in the proof of \Cref{corAnosov}.}
    \label{fig:pseudo_anosov_proof_2}
\end{figure}
 
 Then by the $\delta$-slim triangle condition applied twice, the stability of quasi-geodesics and \Cref{lemma:nearest_point_proj}, we find a point $y$ on one of the sides containing $\alpha'$ that has distance at most $K = C'(\delta,\alpha)+2\delta$ from $f^{N+2M}(\alpha)$ where $C'(\delta,\alpha)$ is the uniform constant coming from quasi-geodesic stability that in turn depends only on $\delta$ and $\alpha$. 
 Moreover, we have (using that $M$ is large) that $y$ splits the geodesic it lies on into subarcs each of length at least $2K$, and hence of length at most $D - 2K$. This is also represented in \Cref{fig:pseudo_anosov_proof_2}.
Then the triangle inequality implies that
 $$d(\alpha', f^{N+2M}(\alpha)) < d(\alpha', f^{N+M}(\alpha)).$$
Now let $f^{N'}(\alpha)$ be the nearest point projection of $\alpha'$ to the quasi-geodesic $ \beta= (f^n(\alpha))_{n \in \mathbb{Z}}$ so that $N' < N$ by the assumption that $\alpha'$ does not belong to $W^N_+$. %\textcolor{red}{shouldn't this be like $N'\leq N+C$ for some $C$? Isn't it true that the preimage of $\{f^n(\alpha)\}_{n\geq R}$ under the nearest point projection map is something like $W_+^{R-C}$ for some $C$?}
We argue now as above. First apply the $\delta$-slim triangle condition to the geodesic triangle with vertices $f^{N+2M}(\alpha), f^{N'}(\alpha)$ and $\alpha'$, as well as the stability of quasi-geodesics to find a point $y$ on the geodesic segments $[f^{N+2M}(\alpha), \alpha']$ or $[f^{N'}(\alpha), \alpha']$ that is distance at most $K' = C'(\delta,\alpha) + \delta$ to the point $f^{N+M}(\alpha)$ (see \Cref{fig:pseudo_anosov_proof_3}). Then using that $M$ is large, so that the distance from $f^{N+M}(\alpha)$ to $f^{N+2M}(\alpha)$ resp.\ $f^{N'}(\alpha)$ is larger than $2K'$, we deduce that $y$ divides the geodesic segment it lies on into segments of length at least $2K'$, and hence of length at most $D-2K'$. Thus we find a (piecewise geodesic) path joining $\alpha'$ to $f^{N+M}(\alpha)$ of length strictly less than $$D' = d(\alpha', f^{N'}(\alpha)) \le d(\alpha', f^{N+2M}(\alpha)) < d(\alpha', f^{N+M}(\alpha)).$$
This is a contradiction to the fact that $f^{N'}(\alpha)$ was a nearest point projection and we conclude that $\alpha' \in W_+^N$.

 \begin{figure}[t]
    \centering
    \includegraphics[width=0.7 \linewidth]{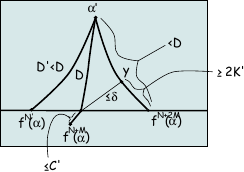}
    \caption{The triangle $\{\alpha',f^{N'}(\alpha),f^{N+2M}(\alpha)\}$. (Here, recall that $K'=C'+\delta$.)}
    \label{fig:pseudo_anosov_proof_3}
\end{figure}
 
 A similar argument (using $g_k^{-1} = \phi^{-1}f^{-k}$ instead of $g_k$ and, for any $j\in\ZZ$,  $\phi^{-1}\circ f^{-j} (\alpha)$ and $\phi^{-1}(W_-^j)$ instead of $f^j(\alpha)$ and $W_+^j$ respectively) shows that any fixed point must also satisfy $\alpha' \in \phi^{-1}(W_-^N)$.
 
 This however contradicts the fact that the two sets $ W_+^N$ and $\phi^{-1}(W_-^N)$ are disjoint. 
 Thus no power of $g_k$ can have fixed points in $X = \mathcal{C}(\Sigma)$ and, in view of \Cref{thm_hyp_Anosov}, we deduce that $g_k$ is pseudo-Anosov for all sufficiently large $k$.
\end{proof}

\begin{remark}\label{rem_Quaaimorphism}
One could also deduce \Cref{corAnosov} as a consequence of results of Bestvina-Fujiwara cf.\ \cite[Theorem 1.1]{Best_Fuji}. In particular, for non-sporadic surfaces, they constructed unbounded quasi-morphisms $h$ on $MCG(\Sigma,\partial \Sigma)$ that are bounded on the stabilizer of any essential simple closed curve. Moreover, one can assume that $h$ is {\bf homogeneous} under (positive) powers, so that $h$ then vanishes on all reducible classes. Consider any $f$ so that $h(f) \ne 0$. Then given any reducible class $\phi$, it follows from the quasi-morphism property that $h(f^k \phi) \ne 0$ for any sufficiently large $k$  and hence $f^k \phi$ cannot be reducible or of finite order by the properties of $h$.  
\end{remark}

\subsection{\textbf{Some hyperbolic geometry}}
\label{sec:hyp_geom}

We recall the following theorem on hyperbolic mapping tori due to Thurston \cite{Th}:
\begin{theorem}[Thurston \cite{Th}]\label{hypmaptori} Let $\Sigma$ be a compact, orientable surface with boundary and negative Euler characteristic. If $\phi$ is a pseudo-Anosov map on $\Sigma$, then the interior of the associated mapping torus has a complete hyperbolic structure of finite volume.
\end{theorem}

We will also need another result due to Thurston on Dehn fillings of hyperbolic manifolds;
an introductory account, as well as a detailed proof, can be found for instance in \cite[Chapter 15]{Mar}.
For the readers' ease, we give here a statement of such a theorem which is adapted to the specific setting in which we will apply it.

Let $N$ be an orientable $3$-manifold with boundary $\partial N$ a finite union $T_1 \sqcup \dots \sqcup T_c$ of $2$-dimensional tori. 
For each $i=1,\ldots,c$, let also $m_i,l_i$ be generators of $\pi_1(T_i)$.
For any $c$-tuple $s=(s_1,\ldots,s_c)$ of \emph{Dehn filling parameters}, i.e.\ of pairs $s_i=(p_i, q_i)$ of coprime integers, one can consider the compact (boundary-less) $3$-manifold $N_{fill}$ obtained by \emph{Dehn filling} the boundary tori with parameters $s=(s_1,\ldots,s_c)$; more explicitly, for each $i=1,\ldots,c$, a solid torus $P_i\coeq\DD^2\times\SSS^1$ is glued to $N$ via the (unique up to isotopy) gluing map $\partial P_i \rightarrow T_i$ sending a meridian of $\partial P_i$ to a curve in the class $p_i m_i +q_i l_i\in\pi_1(T_i)$.

\begin{theorem}[Thurston \cite{ThuNotes}]
\label{HypDehn} 
	In the setting described above, suppose moreover that the interior of $N$ admits a complete hyperbolic metric of finite volume. 
	Then, there is a compact set $K\subset \RR^2$ such that, if every Dehn filling parameter $s_i$ is in $\RR^2\setminus K$, the closed $3$-manifold $N_{fill}$ obtained by Dehn filling $N$ with parameters $s=(s_1,\ldots,s_c)$ admits a finite--volume complete hyperbolic structure $g$. 
	Moreover, the cores of the filling solid tori are closed geodesics of $(N_{fill},g_{hyp})$.
\end{theorem}
Notice that since each $s_i$ is a pair of coprime integers, the theorem implies that one can ensure that a Dehn filling is hyperbolic by excluding finitely many values for {\em each} slope $s_i$.

\begin{remark}\label{inforder} Since the fundamental group of a closed hyperbolic manifold is torsion-free and its closed geodesics are all non-contractible, the cores of the Dehn filling tori will have infinite order in $\pi_1(N_{fill})$. 
\end{remark}

\subsection{\textbf{Proof of the Factorization Lemma}}
\label{sec:fact_lemma}

\begin{proof}(\emph{Factorization Lemma}) Let $f$ be a pseudo Anosov map on $\Sigma$ as in \Cref{corAnosov}. According to \Cref{corAnosov}, $f^k\phi$ is pseudo Anosov on $\Sigma$ for sufficiently large $k$. 
We then write $\phi = F \circ G$, where
$$
F=f^{-k},\; G=f^{k} \phi,
$$ where both are pseudo Anosov for $k\gg 0$. 
By \Cref{hypmaptori}, the interiors of the mapping tori associated to $(\Sigma, F)$ and $(\Sigma, G)$ carry complete hyperbolic structures. 

Let $\gamma_1,\ldots,\gamma_n$ be the components of the boundary $\partial \Sigma$. For each $i=1,\ldots,n$, we then denote by $c_i$ a curve in $\Sigma$ which is parallel to $\gamma_i$ and contained in the interior $\intW$; we can assume, up to isotopy, that they are pairwise disjoint. 
We also denote by $\tau_1,\dots,\tau_n$ the corresponding right-handed Dehn twists, and $\tau\coeq \tau_1\dots \tau_n$.

Observe that $\tau^r=\tau_1^r\dots \tau_n^r$ for every $r \in \mathbb{Z}$, since the $c_i$'s are disjoint.

Let $\phi_1\coeq F\tau^r$ and $\phi_2\coeq \tau^{-r} G$. It is easy to check that the $3$-manifolds $OBD(\Sigma,\phi_1)$ and $OBD(\Sigma,\phi_2)$ correspond to Dehn fillings of, respectively, the mapping tori $\Sigma_F$ and $\Sigma_G$ with respect to Dehn filling parameters $s(r)=(s_1(r),\ldots,s_n(r))$ and $t(r)=(t_1(r),\ldots,t_n(r))$ such that $\norm{s_i(r)},\norm{t_i(r)}\rightarrow +\infty$ for each $i=1,\ldots,n$ as $r\rightarrow +\infty$. 
Thus, for sufficiently large $r$, the hyperbolic Dehn filling Theorem \ref{HypDehn} implies that $OBD(\Sigma,\phi_1)$ and $OBD(\Sigma,\phi_2)$ carry hyperbolic structures and that the binding components (which coincide with the cores of the Dehn filling tori) are geodesics. In particular the latter have infinite order in the fundamental group (see \Cref{inforder}). In other words, we have found the desired decomposition $\phi = \phi_1 \circ \phi_2$ as posited in \Cref{Factlemma}.
\end{proof}

\section{\textbf{Proof of tightness in dimension $5$}}
\label{pfmainthm}

The aim of this section is to prove \Cref{thm:tight_bourgeois} on tightness of the Bourgeois contact structures in dimension $5$. 
For this, we use the following lemma, which is an analogue of the well-known fact that the convex end of a Liouville cobordism with hypertight concave end must be tight \cite{Hof93,AH09}:
\begin{lemma}
	\label{lemma:existence_plane}
	Suppose the connected components of the bindings of $\OBD(\Sigma^2,\phi)$ and $\OBD(\Sigma^2,\psi)$ have infinite order in the corresponding fundamental groups.
	Then, $\BO(\Sigma,\phi\circ\psi)$ is tight.
\end{lemma}

\begin{proof}
	Let $(C,\omega_C)$ be a symplectic cobordism as in \Cref{thm:strong_cobord_bourgeois}. 
	According to \Cref{hypertightobs}, \Cref{item:local_liouv_2} of \Cref{thm:strong_cobord_bourgeois} and our hypothesis on $\OBD(\Sigma,\phi)$ and $\OBD(\Sigma,\psi)$, the Reeb flow of $\lambda_-\vert_{\partial C_-}$ has no contractible periodic orbits. 
	We now show that this implies that $\BO(\Sigma,\phi\circ\psi)$ is tight.
	
	We assume by contradiction that its convex boundary $BO(\Sigma,\phi\circ\psi)$ is overtwisted. 
	According to \cite{BEM}, this implies the existence of an embedded Plastikstufe $\PS$, as defined in \cite{Nie06}.
	Up to attaching a topologically trivial Liouville cobordism to $(C,\omega_C)$ along its positive end, we may then assume that the induced contact form at the positive end is (a positive multiple of) a contact form $\alpha_{PS}$ which is ``adapted'' to $\PS$, i.e.\ it has the normal form described in \cite[Proposition 4]{Nie06} near its core.
	
	Attaching a cobordism at the negative ends using the local Liouville vector fields associated to $\lambda_-$, we obtain the negative Liouville completion $\widehat{C}$ of $\lambda_-$, with a symplectic form $\widehat{\omega}_C$ which coincides with $d (e^t \lambda_-)$ at the negative ends. We now apply the following standard argument. Take an $\widehat{\omega}_C$-compatible almost complex structure $J$, extending the local model of \cite{Nie06}, and cylindrical in the cylindrical ends. 
	We have a Bishop family of Fredholm regular $J$-holomorphic discs in $\widehat{C}$ with totally real boundary, stemming from the core of the Plastikstufe. 
	Analogously to \cite[Proposition 10]{Nie06}, one can check that the exactness of the symplectic form near the positive end, and hence near the Plastikstufe, provides uniform bounds on the Hofer energy, defined as in \cite[Page 115]{WenSFT}. 
	Provided that we rule out boundary bubbling, sphere bubbling and the appearance of symplectic caps, by standard bubbling analysis as in \cite{Hof93,AH09} we can then obtain a finite energy plane in the negative ends. 
	Note that there is no boundary bubbling, as shown in \cite{Nie06}, nor sphere bubbling, by exactness of the symplectic form. 
	We shall need to rule out holomorphic caps, and this can be done as follows (notice that this is not automatic from standard arguments, since $(\widehat{C}, \widehat\omega_C)$ is only pseudo-Liouville). 
	
	Assume the existence of a $J$-holomorphic cap $c$, considered as a map from $\mathbb C$ to $\widehat{C}$. The Hofer energy bounds on the Bishop family provide a periodic Reeb orbit of the Reeb flow of $\lambda_-$, to which $c$ is negatively asymptotic. Now as $\gamma$ is nullhomotopic in $\widehat{C}$, projecting to $\mathbb{T}^2$ via the globally defined projection, we see that the image of $\gamma$ in $\mathbb{T}^2$ is also nullhomotopic. We conclude by \Cref{thm:strong_cobord_bourgeois} that the Reeb orbit $\gamma$ must be a binding component $B_q$. 

	 We now argue that such a cap cannot exist. To this end we let $\pi_\delta \co \widehat C \to C$ be a (smooth) map that collapses the ends of the completed cobordism onto the boundary, depends only on the $t$-coordinate (monotonically) and is the identity on $C$ away from a small $\delta$-neighbourhood of $\partial C$. 
	
	Following \cite[Section 6]{BEHWZ03} we consider the $\omega$-energy of a holomorphic plane
	$$
	\mathbf E_\omega(c)= \int_\mathbb C c^*\overline \omega_\delta,
	$$
	where $\overline \omega_\delta = \pi^*_\delta\omega_C$ is the $2$-form which is explicitly given by
	$$
	\overline\omega_\delta=\left\{\begin{array}{cc}
	 d\lambda_+    &  \mbox{ on } [0,\infty)\times \partial_+ C\\
	 \pi_\delta^*\omega_C    & \mbox{ on $C$} \\
	 d\lambda_- & \mbox{ on } (-\infty,0]\times \partial_- C.
	\end{array}\right.
	$$
	In fact the resulting forms are all cohomologous relative to the ends for any $\delta>0$ so that the precise value is irrelevant. Taking a limit as $\delta\searrow 0$, we obtain the (piecewise) description of the $\omega$-energy of \cite{BEHWZ03}.
	
	As the almost complex structure is cylindrical on the ends, it is easy to check that $\mathbf{E}_{\omega}(c)$ is non-negative, and strictly positive if $c$ is not completely contained in an end. In particular, given the assumption that all Reeb orbits have infinite order in the fundamental group of the negative end, we deduce that the holomorphic cap cannot lie completely in the negative cylindrical end, and hence $\mathbf E_{\omega}(c) > 0$.
	
	By \Cref{item:sympl_form_2} of \Cref{thm:strong_cobord_bourgeois}, there is a primitive $\nu$ of $\omega_C$ on the initial cobordism that is positive along binding components. We remark that this primitive will not be of the form $d(e^t \lambda_-)$ near the negative end, as the cobordism is only pseudo-Liouville. Pulling back under the map $\pi_\delta$ that collapses each end in turn gives a primitive for $\overline \omega_\delta$ that naturally extends to the compactification given by adding $\{\infty\} \times M_{\pm}$, since it is invariant under positive translations in the $t$-direction. We denote the resulting primitive for $\overline \omega_\delta$ by $\overline \nu$.
	
	Then integrating the exact form $\overline \omega_\delta= d \, \overline \nu$ along the holomorphic cap $c$, and using Stokes' theorem, we obtain: 
	\begin{equation*}
	0<\int_{c}\overline \omega_\delta = 
	\int_{-\gamma}\overline \nu <0
	\end{equation*}
	and this contradiction finishes the proof.
	%We conclude that the $J$-holomorphic plane obtained above is contained in the bottom level, and therefore we obtain a contractible Reeb orbit in the negative symplectization of $\lambda_-$. But there are no such orbits at the negative ends by our assumptions, .
\end{proof}

We can now proceed to the proof of \Cref{thm:tight_bourgeois}:

\begin{proof}[Proof (\Cref{thm:tight_bourgeois})]
	
	We start by proving tightness in the ``generic'' case of non-sporadic page $\Sigma^2$. We then deal with sporadic pages on a case by case basis. 
	Lastly, we explain how to deduce universal tightness by the proof of tightness.
	\medskip
	
	\textbf{Case 1: non-sporadic $\Sigma$.} By \Cref{Factlemma} we may factorise the monodromy $\phi=\phi_1 \circ \phi_2$, where the components of the bindings in $OBD(\Sigma,\phi_1)$ and $OBD(\Sigma,\phi_2)$ have infinite order. Then, according to \Cref{lemma:existence_plane}, we conclude that $BO(\Sigma,\phi)$ is tight.
	
	\medskip
	
	\textbf{Case 2: $\Sigma$ is a disc.} In this case, the monodromy $\phi$ is necessarily isotopic to the identity. 
	In other words, the resulting contact $3$-manifold is $(\SSS^3,\xistd)$ and the open book structure is the one induced by the subcritical Stein-filling $\DD^4$. 
	According to \cite[Theorem A.(b)]{LMN}, the associated Bourgeois contact structure is Stein fillable, and hence tight.
	
	\medskip
	
	\textbf{Case 3: $\Sigma$ is an annulus.} The mapping class group of the annulus is generated by a single positive Dehn twist around the core circle. 
	If the monodromy is a non-negative power of such generator, then the resulting contact $3$-manifold is Stein fillable;
	then, according to \cite[Example 1.1]{MNW}, the associated Bourgeois contact structure is weakly fillable, and hence tight. 
	If the power is negative, according to \cite[Theorem B]{LMN}, the Bourgeois contact structure associated to $OBD(\Sigma,\phi)$ is contactomorphic to that associated to $OBD(\Sigma,\phi^{-1})$, so we obtain tightness for this case. 
	
	\medskip
	
	\textbf{Case 4: $\Sigma$ is a pair of pants.}
	For simplicity, enumerate from $1$ to $3$ the connected components of $\partial \Sigma$. 
	For $i=1,2,3$, let $\tau_i$ be a positive Dehn twist along the $i$-th connected component of $\partial \Sigma$; these give generators for the abelian mapping class group of the pair of pants (cf.\ \cite[Section 3.6.4]{FM}). Thus we can write $\phi =\tau_1^{a_1}\circ\tau_2^{a_2}\circ \tau_3^{a_3}$.
	We then set $\tau\coeq \tau_1\circ\tau_2\circ\tau_3$ and, for any $N\in\NN_{>0}$, we can decompose $\phi$ as $\phi = F \circ G$, with $F\coeq \phi \circ \tau^N= \prod_{i=1}^3 \tau_{i}^{N+a_i}$ and $G\coeq \tau^{-N}=\prod_{i=1}^{3}\tau_i^{-N}$.
	We then use the following result:
	\begin{lemma}
		\label{lemma:pair_pants}
		If $N>0$ is big enough, each binding component of $\OBD(\Sigma,F)$ is of infinite order in $\pi_1\left(\OBD\left(\Sigma,F\right)\right)$.
		The same is true for $\OBD(\Sigma,G)$.
	\end{lemma}
	
	\noindent
	Using \Cref{lemma:existence_plane,lemma:pair_pants}, we conclude that $\BO(\Sigma,\phi)$ is tight, as desired.
    It's then only left to prove \Cref{lemma:pair_pants} in order to conclude the proof of tightness in the case of a pair of pants:
\begin{proof}[Proof (\Cref{lemma:pair_pants})]
	We deal only with the case of $\OBD(\Sigma,F)$; the proof for the manifold $\OBD(\Sigma,G)$ is completely analogous.
	
	We first point out that, as explained in detail for instance in \cite[Section 3]{Ozb07},  the manifold $\OBD(\Sigma,F)$ can be seen as obtained by Dehn surgery on the total space of $\SSS^2\times\SSS^1\rightarrow\SSS^2$ along three $\SSS^1$-fibers, with coefficients $r_i\coeq-\frac{1}{N+a_i}$, for each $i=1,2,3$. In other words, $\OBD(\Sigma,F)$ is the Seifert manifold 
	\begin{equation*}
	\{0,(o_1,0);(N+a_1,-1),(N+a_2,-1),(N+a_3,-1)\} \text{ .}    
	\end{equation*}
	Moreover, the orbit space $O$ of the Seifert fibration of $\OBD(\Sigma,F)$ is a $2$-dimensional orbifold, with underlying topological surface $\SSS^2$, and the binding $B$ of $\OBD(\Sigma,F)$ consists of a union of fibers of the Seifert fibration.
	
	We recall that there is a notion of \emph{orbifold Euler characteristic} $\chi_{orb}$ for orbifolds that behaves multiplicatively under finite covers of orbifolds. 	In our special case of the base orbifold $B$ of the Seifert fibered space $\OBD(\Sigma,F)$, we have
	\begin{equation*}
	\chi_{orb}(O)=\chi(\SSS^2) - \sum_{i=1}^{3}\left(1+\frac{1}{N+a_{i}}\right)=-1-\sum_{i=1}^{3}\frac{1}{N+a_{i}} \text{ .}
	\end{equation*}
	From now on, let $N>0$ be so big that $\chi_{orb}(O)<0$.
	In particular, $\OBD(\Sigma,F)$ is finitely covered by a circle bundle $X$ over a hyperbolic surface $S$, in such a way that fibers of $X\rightarrow S$ are mapped to fibers of $\OBD(\Sigma,F)\rightarrow O$ (see \cite{Sco83} for instance). 
	Now, $S$ being hyperbolic, the fibers of $X$ are of infinite order in $\pi_1(X)$. 
	As $X$ covers $\OBD(\Sigma,F)$ in a compatible way with their Seifert bundle structures, it follows the fibers of $\OBD(\Sigma,F)$, hence its binding too, are of infinite order in its fundamental group, as desired.
\end{proof}

\textbf{Universal tightness.}
Note that all the above arguments remain valid for the pull-back under any finite cover of $V \times \TT^2$ induced by a finite cover of the first factor. 
Since finite covers over the second factor do not change the contact structure up to contactomorphism (cf.\ Remark \ref{rem:cover_finite}), such covers also preserve tightness. 
Now {\em any} finite cover is itself covered by a composition of covers of the respective factors.
Consequently, the contact structure remains tight under any finite cover on the first factor. 
Since the fundamental group of any closed 3-manifold is residually finite (cf. \cite{Hempel}) so is $\pi_1 (V \times \TT^2)$ and hence tightness on finite covers is equivalent to tightness on the universal cover of $V \times \TT^2$ and universal tightness follows. 
This concludes the proof of \Cref{thm:tight_bourgeois}. 
\end{proof}

%$\hfill \square$

\section{\textbf{Obstructions to symplectic fillability }}
\label{sec:general_obstr_fillab}
In this section we describe a general capping construction for $\mathbb S^1$-invariant contact structures by applying the handle attachments of Massot-Niederkr\"uger-Wendl \cite{MNW}, and use it to give new obstructions to symplectic fillability in all dimensions.

\subsection{{\bf $\SSS^1$--invariant contact structures in terms of Giroux domains}}\label{sec_S_inv} 
%We want to argue here that the $\SSS^1-$invariant contact structures constructed in \cite[Theorem 1]{GeiSti10} can be reinterpreted in terms of Giroux domains as defined in \cite[Section 5]{MNW}.
%In order to do so, let us start by recalling the following definition. 
%{\color{red} What was written here was wrong!!! (now commented out)}
As discussed in \cite[Section 6]{DinGei12} any $\SSS^1-$invariant contact structure induces a decomposition into so called ``ideal Liouville domains'' as defined in \cite[Section 5]{MNW}. 
We now describe how this decomposition arises. First of all we recall:

\begin{definition}[Giroux \cite{GirouxIdealLiouvDom}]
    An \emph{ideal Liouville domain} $(V,\omega)$ is a domain $V$ with an exact symplectic form $\omega$ defined on the interior of $W$ that admits a primitive $\lambda$, satisfying the following condition: for each (and hence any) function $u\colon V\to\RR_{\geq0}$ with regular level set $\partial V = \{u=0\}$, the form $u\lambda$ extends smoothly to $V$ inducing a contact form on $\partial V$.
\end{definition}
Notice that, by \cite[Proposition 2]{GirouxIdealLiouvDom}, at the boundary of any ideal Liouville domain $(V,\omega)$ there is a contact \emph{structure} $\xi$ which is well defined, i.e.\ only depending on $\omega$.
Moreover, according to \cite[Example 8]{GirouxIdealLiouvDom}, the complement of the dividing set of a convex hypersurface in a contact manifold is a disjoint union of ideal Liouville domains (the orientation of each of which can either correspond or be opposite to that induced by the ambient orientation and the transverse contact vector field on the hypersurface).

Now, given an ideal Liouville domain $(V,\omega)$, one can consider the associated \emph{Giroux domain} as defined in \cite[Section 5.3]{MNW}, namely the contact manifold $(V\times \SSS^1,\ker(u dt + u \lambda))$; this is well defined in the sense that, up to isotopy relative boundary, it only depends on $\omega$ and not on the function $u$ nor the primitive $\lambda$.

We now consider an $\SSS^1-$invariant contact structure on a product $V\times\SSS^1$. Note that the hypersurfaces $V\times \{pt\}$ are convex in the sense of Giroux \cite{Gir91}, i.e.\ transverse to the contact vector field $\partial_\theta$, where $\theta\in\SSS^1$.
Now, according to what we pointed out above, this determines a decomposition $V=V_{+}\cup \overline{V}_-$ with $V_+,V_-$ ideal Liouville domains (where on $V_-$ the ideal symplectic form is negative), where $V_\pm$ is the set where $\partial_\theta$ is positively/negatively tranverse to the contact structure. Moreover, by $\SSS^1-$invariance, the contact structure on $V\times\SSS^1$ is just obtained by gluing the two Giroux domains $V_+\times \SSS^1$ and $\overline{V}_-\times \SSS^1$ associated to $V_+$ and $\overline{V}_-$ respectively along their common boundary $N\times \mathbb S^1$, with $N\coloneqq \partial V_+ = \partial \overline{V}_-$, the dividing set, along which $\partial_\theta$ is tangent to the contact structure.

\medskip
\textbf{A capping cobordism.} 
Let now $W^{2n+2}$ be a hypothetical strong symplectic filling of $V\times \SSS^1$.
The aim of this subsection is to attach a symplectic cobordism on top of $W$ which, near its positive boundary, splits as the product of a symplectic $\SSS^2$ and portion of the symplectization of $\partial V_+$, where $V=V_+\cup\overline{V}_-$ is the decomposition in Liouville domains of the convex hypersurface $\{pt\}\times V$ in the $\SSS^1-$invariant contact manifold $V\times \SSS^1=\partial W$.
In order to achieve this, we use the symplectic cobordism described in \cite[Section 6]{MNW}.

First, we describe what the cobordism should look like topologically (i.e.\ without specifying the symplectic form or dealing with details such as corner smoothings), as also done in \cite[Section 6.1]{MNW} in their analogous setting.
For this, it is more convenient to see $V\times \SSS^1$ as obtained as the union $V_+\times \SSS^1 \cup [-\delta,\delta]\times N \times \SSS^1 \cup \overline{V}_-\times \SSS^1$, as explained above.
Now, the cobordism is obtained by attaching $V_\pm\times \DD^2_\pm$ on top of $V_\pm\times \SSS^1\subset V\times \SSS^1$  (and smoothing corners), where the latter is seen as the positive boundary of the trivial cobordism $[0,1]\times V\times \SSS^1$.
The positive boundary component of the resulting cobordism is then just $N\times \SSS^2$.
Notice also that there are distinguished submanifolds $C_{\pm}\coloneqq V_\pm\times\{0\}\subset V_\pm\times \DD^2$ in this capping cobordism, with the property that, once removed, the cobordism deformation retracts onto the negative boundary.

Now that the topological picture has been described explicitly, we need to argue that the cobordism has a symplectic structure suited to our needs. 
In fact, this is precisely the content of \cite[Section 6]{MNW}, where the smoothing procedure and the symplectic structure are described with great care. We hence refer the reader to that paper for details, and limit ourselves to stating the needed properties of the cobordism in the following lemma.

\begin{lemma}
\label{lemma:capping_cobord_high_dim}
There is a symplectic cobordism $(C,\omega_C)$ with negative (i.e.\ concave) contact boundary $\partial_- C=V\times \SSS^1$ and weakly convex positive boundary
\[
\partial_+ C = N \times \SSS^2
\] 
where $N=\partial V$.
Moreover, there is a tubular neighborhood $(-\delta,0]\times\partial_+ C$ such that $\omega_C$ is of the form $d(e^t\alpha_N)+\omega_{S}$, where $t\in(-\delta,0]$ and 
\begin{itemize}
    \item[$\bullet$] $\omega_S$ is an area form on $\SSS^2$,
    \item[$\bullet$] $\alpha_N$ is a contact form on $N$.
\end{itemize}
\noindent
Lastly, there are symplectic submanifolds $C_\pm$, diffeomorphic to $V_\pm$, such that $C\setminus C_\pm$ deformation retracts onto its negative boundary, and such that $C_\pm$ intersect transversely, positively and in exactly one point, each symplectic sphere in the previously described neighborhood of the positive boundary $\partial_+ C$.
%$(C,\omega_C)$ can be interpreted as attaching a product symplectic handle $H = X\times D_- \sqcup X\times D_+$ on top of the compact piece of symplectization $[0,1]\times BO(\Sigma,\phi)$, where $\phi$ is seen as .
\end{lemma}

\noindent In what follows, we denote by $\Wcap$ the result of stacking $(C,\omega_C)$ on top of the strong symplectic filling $W$ along the common boundary $V\times \SSS^1$.

\medskip
\noindent
\textbf{A moduli space of spheres in the capped filling.}
We now consider a moduli space of pseudo-holomorphic spheres in $\Wcap$.
The setup and properties needed are essentially the same as those discussed in \cite[Section 7.2]{MNW}. For this reason, we limit ourselves to describing the situation and stating the necessary properties here, and refer the reader to \cite{MNW} for further details. Note however that in their situation the positive boundary of the capped filling has two connected components, one of contact convex type and the other of stable Hamiltonian type, which leads to a contradiction; in our setting, we instead end up with one boundary component.

Now according to \Cref{lemma:capping_cobord_high_dim}, the boundary $\partial\Wcap$ has a collar neighborhood in $\Wcap$ so that the symplectic form is split and of the following form
\begin{equation*}
    \left((-\delta,0]\times N \times \SSS^2, d(e^{t}\alpha)  + \omega_{S}\right), \textrm{ where }  t\in(-\delta,0].
\end{equation*}
Hence, one can find an almost complex structure $J$ compatible with the symplectic structure which, on this neighborhood, splits as a direct sum of an almost complex structure  $J_{S}$ on the 2-sphere and some $J_{\alpha}$ on the symplectization factor $((-\epsilon,0]\times N, d(e^t\alpha))$.
In particular, for any point $(t,q)\in (-\epsilon,0]\times N$, the sphere
\begin{align*}
    u_{t,q} \co  \SSS^2 & \to (-\epsilon,0]\times N\times \SSS^2 \\
    y & \mapsto (t,q,y)
\end{align*}
is $J$--holomorphic.
Moreover, as the co-cores $C_\pm$ from \Cref{lemma:capping_cobord_high_dim} are symplectic submanifolds of $\Wcap$, one can choose such a $J$ in such a way that $C_\pm$ are $J$-invariant in $(\Wcap,J)$; we hence assume this is the case.
Lastly, in the complement of a neighborhood of $\partial\Wcap\cup C_+ \cup C_-$, the $J$ can then be chosen to be generic, so that every simple holomorphic sphere intersecting said region is regular. 

We also have regularity near the boundary, which is proved in \cite[Page 334]{MNW} in a completely analogous setting:

\begin{lemma}[Fredholm Regularity]\label{lemma:reglemma}
	The spheres $u_{t,q}$ are Fredholm regular and have Fredholm index $\dim(X) =2n$.
\end{lemma}
%\begin{proof} 
%	Fix $u=u_{t,q}$. Since the almost complex structure is split near the boundary and the projection to $(-\epsilon,0]\times \partial X$ is $J$-holomorphic, the normal bundle $N_u$ to $u$ is the pullback of $T_{(t,q)}((-\epsilon,0]\times \partial X)$ via the projection to $(-\epsilon,0]\times\partial X$. 
%	In particular, a trivialization of this vector space via a complex basis pulls back to a trivialization of the bundle $N_u$, which hence splits as a sum of holomorphic line bundles.
%	This implies that the normal component of the linearized Cauchy-Riemann operator splits as $\mathbf{D}_u^N=\overline{\partial} \oplus \cdots \oplus \overline{\partial}$, where $\overline{\partial}$ is the standard Cauchy-Riemann operator acting on $W^{1,2}$-sections of the corresponding line bundle. 
%	By automatic transversality, each summand is surjective, and we deduce surjectivity of $\mathbf{D}_u^N: W^{1,2}(N_u)\rightarrow L^2(\Omega^{0,1}(N_u))$, so that $u$ is indeed Fredholm regular. 
%	
%	Since $u$ is embedded, its index is the Fredholm index of $\mathbf{D}_u^N$, and from the Riemann-Roch formula we obtain
%	\begin{equation*}
%	\mbox{ind}(u)=\mbox{ind}(\mathbf{D}_u^N)=\mbox{rank}_\mathbb{C}(N_u)\chi(u)+2\underbrace{c_1(N_u)}_{=0} =2n \text{.}
%	\qedhere
%	\end{equation*}
%\end{proof}

One can then consider the connected component of the moduli space $\calM$ of $J$-holomorphic spheres in $\Wcap$ containing the pseudo-holomorphic spheres $u_{t,q}$. 
%By uniqueness near the boundary of $\Wcap'$, the moduli space $\calM$ %extends the family of $u_{t,q}$'s defined on the cylindrical end over $\partial \Wcap'$. 
%is a foliation near $\partial \Wcap'$.
We let $\calM_*$ denote the corresponding marked moduli space on which there is a naturally defined evaluation map $ev \co \calM_* \to \Wcap$.
Lastly, we denote by $\calMbar, \calMbar_*$ the Gromov compactification of $\calM, \calM_*$ respectively; abusing notation slightly, we also continue to denote the evaluation map on the compactification by $ev \co \calMbar_* \to \Wcap$.

Now, we have the following important uniqueness property near the boundary, which can also be proven exactly as in \cite[Pages 334 and 335]{MNW}:
\begin{lemma}[Local Uniqueness]\label{lemma:localuniq} For sufficiently small $\epsilon>0,$ any curve $u$ in the moduli space $\overline{\mathcal{M}}$, that intersects the collar $\mathcal{N=}(-\epsilon,0]\times N \times \SSS^2$, is a reparametrization of one of the spheres $u_{t,q}.$
\end{lemma}

In particular, $\overline{\mathcal{M}}, \overline{\mathcal{M}}_*$ are smooth manifolds near their boundary, and the restriction $ev_\partial$ of $ev\colon \calMstarbar\to \Wcap$ to a (sufficiently small) neighborhood of $\partial \calMstarbar$ is a diffeomorphism onto a neighborhood of $\partial \Wcap$. 
We also have a diffeomorphism $\partial \overline{\mathcal{M}}\cong N$. It follows that each of the curves in $\calM$ has to be simple, as they intersect each co-core geometrically once by positivity of intersections together with uniqueness near the boundary (as the latter implies that intersections cannot escape at the boundary).

\medskip

\noindent \textbf{Nodal stratification: Semi-positive case.} In this case the moduli space $\overline{\mathcal{M}}$ is a stratified space%
\footnote{In this paper, by \emph{stratifed space}, we will mean a filtration $\emptyset=\mathcal{M}^{m+1}\subset\mathcal{M}^m\subset\dots\subset \mathcal{M}^0=\mathcal{M}$ of a compact topological space $\mathcal{M}$, where the interior $\mbox{int}(\mathcal{M}^i)$ of each $\mathcal{M}^i$ is called the \emph{$i$-th stratum}, and the closure of each stratum satisfies $\overline{\mathcal{M}}^i= \bigcup_{j\geq i}\mathcal{M}^j$.} (cf.\ \cite[Section 6.5]{McDSal}):
$$
\emptyset=\overline{\mathcal{M}}^{n+1}\subset \overline{\mathcal{M}}^{n}\subset \cdots  \subset  \overline{\mathcal{M}}^1 \subset \overline{\mathcal{M}}^0= \overline{\mathcal{M}}
.$$
Here, $\overline{\mathcal{M}}^i$ consists of nodal configurations of spheres which have at least $i$ nodes, and the interior $\mbox{int}(\mathcal{M}^i)$, consists of nodal configurations with precisely $i$ nodes. 
For our choice of $J$, the top open stratum $\mathcal{M}=\mbox{int}(\overline{\mathcal{M}}^0)$ is a smooth $2n$-dimensional manifold consisting of simply covered spheres. The elements of the $i$-th stratum int$(\overline{\mathcal{M}}^i)$ contain \emph{main} sphere components, which intersect at least one of the $J$-invariant co-cores $C_\pm$ in at most once and so are again simply covered, while the rest are possibly multiply covered bubbles which intersect no co-cores. 
By the Uniqueness \Cref{lemma:localuniq}, no element in $\overline{\mathcal{M}}^i$ touches the boundary of $W_{cap}$, for $i\geq 1$. 

Similar remarks hold for the marked moduli space, which is a stratified space $$\emptyset=\overline{\mathcal{M}}_*^{n+1}\subset \overline{\mathcal{M}}_*^{n} \subset \cdots \subset  \overline{\mathcal{M}}_*^1 \subset \overline{\mathcal{M}}_*^0= \overline{\mathcal{M}}_*,$$ and the forgetful map respects the stratification.

The semi-positive assumption implies that the dimension of the image of $\overline{\mathcal{M}}^i_*$ under the evaluation map is \emph{at most} $2n+2-2i$, in the sense of \cite[Sec.\ 6.5]{MS}. More precisely, the image under the evaluation map of $\overline{\mathcal{M}}^i_*$ is covered by the images of underlying moduli spaces of \emph{simple} stable maps, each of which is actually a smooth manifold for generic choice of almost complex structure. 
With a slight abuse of language, we shall occasionally say that the corresponding unmarked piece $\overline{\mathcal{M}}^i$ has dimension at most $2n-2i$.

\begin{remark}\label{rem:semi_pos_6_dim}
Since the semipositivity condition is automatically satisfied for every $6$-dimensional symplectic manifold (see e.g.\ \cite[Section 6.4]{McDSal}), in this case the evaluation map on the marked moduli space determines a pseudo-cycle on $\Wcap$ (see e.g.\ \cite[Section 6.5]{McDSal}).
\end{remark}

%\begin{remark}\label{rem:multdisj} \textcolor{red}{Is this remark anymore relevant, if we deal with semipositive case instead of sympl aspherical?} In the general case where $W$ is an arbitrary strong filling, there might be elements in $\calMbar\backslash \calM$ with multiply covered sphere components. However, since every element in $\calMbar$ intersects the co-cores once and positively, those multiply covered components are necessarily disjoint from the co-cores.
%\end{remark}

\medskip\noindent
\subsection{\textbf{Proof of \Cref{thm:homology_injection_boundary}}}

We are now ready to give a proof of \Cref{thm:homology_injection_boundary} from the Introduction, i.e.\ of the fact that the inclusion of $N = \partial V_\pm $ induces an injection in rational homology for classes that survive in $V_\pm$.
We prove this first in the semi-positive setting, which uses the more standard language of pseudo-cycles as documented for instance in \cite[Section 6.5]{McDSal}.
%(Recall however that in the setting of dimensions $\leq6$, i.e.\ $n=2$ for us, semi-positivity automatically holds.)
Then, we prove the result in the general case using the polyfold technology developed by Hofer--Wysocki--Zehnder in the Gromov-Witten case in \cite{HWZ_GW}.
\\

%We begin with the case that the filling $W$ is semi-positive:

\begin{proof}[Proof of \Cref{thm:homology_injection_boundary} in the semi-positive case]
We start by fixing some notation. 
Let $i'$ be inclusion
\[
N = N\times\{pt\}\hookrightarrow \partial \Wcap = N\times \SSS^2
\]
for $pt$ a point in the equator of $\SSS^2$.
%$$\Sigma\hookrightarrow \partial X = \mathcal{D}_{\Id}\Sigma \times\SSS^1_x \overset{i''}{\hookrightarrow} \partial \Wcap' = \partial X \times \SSS^2,$$ 
%where the first arrow is the inclusion of $\Sigma$ as one of the two halves of the double $\mathcal{D}_{\Id}\Sigma$ and $i''$ is the inclusion of $\partial X$ into $\partial X \times \SSS^2$ as $\partial X \times \{pt\}$ with $pt$ a point in the equator of $\SSS^2$.
We also let $i$ denote the composition of $i'$ with $\partial \Wcap \hookrightarrow \Wcap$.
(Notice that $i(N)$ in fact lives in $\partial W\subset \Wcap$.)

By \Cref{lemma:capping_cobord_high_dim}, we have that each element in the moduli space $\calMbar$ intersects each co-core $C_{\pm}$ transversely in a single point. 
Taking the intersection of a curve with the co-core $C_\pm$ then gives a continuous map $\mathcal{I}^\pm: \calMbar\to C_\pm$, which, by uniqueness of the $u_{t,q}$'s near the boundary $N \times \SSS^2$ of $\Wcap$, is a diffeomorphism near $\partial C_\pm$. 

Define now the maps $\mathcal{I}_*^\pm$ and $f_\pm$ so that the following diagram is commutative:
\begin{equation}
\label{eqn:diagram_pf_homol_inject}
	\begin{tikzcd}
		N \arrow[r, hook, "i'"] \arrow[rrrr, "i", bend left=23] \arrow[rrrdd, hook, "f_\pm"] & 
		N \times \SSS^2 = \partial \Wcap \arrow[r,"ev_\partial^{-1}", "\cong"'] &
		\partial \calMstarbar \arrow[r, hook, "j"] &
		\calMstarbar \arrow[r, "ev"] \arrow[d, "\pi"'] \arrow[dd, bend left=30, "\calI_*^\pm"] &
		\Wcap \\
		&
		&
		&
		\calMbar \arrow[d, "\calI^\pm"'] & 
		\\
		&
		&
		&
		C_\pm \simeq V_\pm & 
	\end{tikzcd}    
\end{equation}
Here, $j$ is the 
natural inclusion $\partial \calMstarbar \subset \calMstarbar$ and $\pi$ is the map forgetting the marked point.
(Recall also that $ev$ is a diffeomorphism near the boundary, by the local uniqueness \Cref{lemma:localuniq}.)
Moreover, by the explicit definitions of $i'$ and $\mathcal{I}^\pm$, the map $f_\pm$ is actually homotopic to the natural inclusion of $N=\partial V_\pm$ into $C_\pm\simeq V_\pm$.
We lastly note that, because we are in the semi-positive case, the evaluation map gives a pseudo-cycle representing a relative fundamental class that we denote $[\overline{\mathcal{M}}_*]$. 

Recall that the statement of \Cref{thm:homology_injection_boundary} claims that the natural composition $V_\pm \simeq V_\pm \times\{pt\}\hookrightarrow \partial W = V\times\SSS^1\hookrightarrow W$ induces an injective map on rational homology once restricted to $\image((I_\pm)_*)$, where $I_\pm\colon N=\partial V_\pm\hookrightarrow V_\pm$ is the natural injection, and $(I_\pm)_*$ is the corresponding map in rational homology.
In order to prove this, it is enough (and in fact stronger) to prove in fact that the same conclusion holds for the composition $V_\pm\hookrightarrow W \hookrightarrow \Wcap$.
We hence aim at proving the property for this last map then.

%We now point out that the image of the inclusion $i$ described above is disjoint from the co-cores $C_{\pm}$. 
% (recall \Cref{fig:handleS1eq2}). 
%Moreover, 
Notice that one can find a deformation retraction of $\Wcap$ minus slight push-offs of $C_{\pm}$ onto $W$ which pushes the co-cores $C_\pm$ to $V_\pm\times\{pt\}\subset V\times \SSS^1 = \partial W$.
In particular, 
%one can see that the composition $i$ of $i'$ with
%$\partial \Wcap\hookrightarrow \Wcap$ is homotopic to the composition of the inclusion $N=\partial V_\pm \hookrightarrow V_\pm \hookrightarrow W$ described in the statement of \Cref{thm:homology_injection_page} with the natural inclusion $W\hookrightarrow\Wcap$.
%In particular, 
in order to prove statement of \Cref{thm:homology_injection_boundary} %(for $V_\pm\hookrightarrow W \hookrightarrow \Wcap$, by what said above)
it is enough to prove the following: given a class $z\in H_*(N;\QQ)$ such that $(f_\pm)_*(z)\neq 0$ in $H_*(C_\pm;\QQ)$, then $i_*(z)\neq 0$ in $H_*(\Wcap;\QQ)$.
By looking at its contrapositive, we hence need to prove that, if $i_*(z)=0$ in $H_*(\Wcap;\QQ)$, then $(f_\pm)_*(z)= 0$ in $H_*(C_\pm;\QQ)$ as well.

%$f$ in the above diagram is homotopic, as a map with values in 
%in order to prove \Cref{thm:homology_injection_boundary}, it is enough to prove the following:
%let $z$ be a non-trivial homology class in 
%that $i$ is injective in rational homology. 
%We lastly note that, because we are in the semi-positive case, the evaluation map gives a pseudo-cycle representing a relative fundamental class that we denote $[\overline{\mathcal{M}}_*]$. 

%We now prove that $i$ is injective in rational homology by contradiction.
Let then $(f_\pm)_*[\sigma]\in H_*(C_\pm;\QQ)$ be non-trivial, where $\sigma$ is rational cycle in $N$, such that $i_*[\sigma]=0$ in $H_*(\Wcap;\QQ)$.
We then write $\sigma=\partial b$ as a boundary in $\Wcap$.

In the ideal case that all spaces and cycles are actually smooth manifolds one could simply finish the proof as follows. 
After a small perturbation to ensure transversality, the preimage $ev^{-1}(b)$ in $\overline{\calM}_*$ is a chain with boundary $ev^{-1}(\sigma)$. 
As $ev_\partial$ is a diffeomorphism and $\sigma$ (seen in $\Wcap$) lies in its image, this implies that $(\mathcal{I}_\pm)_*(ev^{-1}(b))$ has boundary $\sigma$ in $C_+$, thus proving that $(f_\pm)_*[\sigma]=0$ in $H_*(C_\pm;\QQ)$.
However, since $\overline{\calM}_*$ is not in general a manifold, we need rely on pseudo--cycles, as follows.

We first claim that the map $(ev,\calI_*)\co \calM_* \to \Wcap \times C_+'$ is a $(2n+2)$--dimensional pseudo--cycle. 
	Indeed, the closure of the image of $(ev,\calI_*)$ is the image of the compactified moduli space $\overline{\calM}_*$, on which the map $(ev,\calI_*)$ is continuous. 
	Moreover, the image of the nodal set $\overline{\mathcal{M}}_*^1$ for both maps agrees with the image of the space $\overline{\mathcal{M}}_*^s$ of underlying simple curves. 
	As $\overline{\mathcal{M}}_*^s$ consists of components of dimension $\leq 2n$ by semi-positivity, this gives the claim.
	
    Now, one can perturb $b$ (relative to its boundary) so that $b \times C_+$ becomes weakly transverse to the pseudo--cycle $(ev,\calI_*)$, in the sense of \cite[Definition 6.5.10]{MS}. 
    Strictly speaking, in order to apply their definition to our setting, the cycles $b$ and $\sigma$, would need to be represented by embedded submanifolds (up to integer multiples). This will however only be true for (relative) cycles in degree $< n$. 
    In the general case they can be represented, according to a result of Thom \cite{Tho54}, by smooth maps from manifolds, possibly with boundary, to the target space.
    This being said, all the considerations in \cite[Section 6.5]{McDSal}, including the definition of weakly transverse and the proposition that we use below, carry on with the same proof to this more general setting.
    To ease the notation, we will however omit the explicit mention of these maps representing $\sigma,b$.
    
	According to (a relative version of) \cite[Proposition 6.5.17]{MS}, the restriction of $\calI_*$ to $ev^{-1}(b)$ gives a relative pseudo--cycle on $C_+$. 
	(Here $ev^{-1}(b)$ should be interpreted as the preimage of the diagonal under the map $(b,ev):B \times \calM_* \to W \times W$ in the case the the representing map $b\co B \to W$ is not a submanifold.)
	%(Here the pre-image $ev^{-1}(b)$ is to be interpreted as the sum of the simplex-wise pre-images, which then glue to give a relative pseudo-cycle; and similarly for $ev^{-1}(\sigma)$, which is the boundary of $ev^{-1}(b)$.)
	Moreover, the relative pseudo--cycle $\calI_*(ev^{-1}(b))$ is a weak representative of a relative homology class (see \cite[Discussion after Lemma 6.5.6]{MS}). 
	Finally, as $ev$ is a diffeomorphism near the boundary, the boundary of the relative pseudo-cycle $\calI_*(ev^{-1}(b))$ in $C_+$ is the image of the original homology class $[\sigma]$ via $(f_\pm)_*\colon H_*(\Sigma;\QQ)\to H_*(C_+;\QQ)$. 
	In particular, $(f_\pm)_*[\sigma]=0$ in $H_*(C_\pm;\QQ)$, as desired.
	
	The statement about surjectivity at the level of the fundamental group follows, as pointed out to us by a referee, from the fact that bubbling occurs in codimension 2 for semi-positive fillings, again via a push-pull argument. 
	This can be carried out by an argument as in the proof of Lemma 2.3 in  Ghiggini-Niederkrüger-Wendl \cite{GNW}, to which we refer the reader for more details.	
	%$f_\pm\colon \Sigma\to C_+$ is not injective in rational homology, contradicting the fact that it is homotopy equivalent to the inclusion $\Sigma\hookrightarrow C_+=\Sigma \times D^*\SSS^1_x$.
	\end{proof}

We now need to deal with the general setting of possibly non semi-positive fillings. As in the proof in the semi-positive case above, we have that each element in the moduli space $\calMbar$ intersects each co-core $C_{\pm}$ transversely in a single point. 
This again gives a map $\mathcal{I}: \calMbar\to C_+$, which, by uniqueness of the $u_{t,q}$'s, is a diffeomorphism near the boundary. 
However, as we are in the higher dimensional setting, the evaluation map on the marked and compactified moduli space will not automatically be a pseudo-cycle (e.g.\ in dimensions $2n+2\geq 8$). 
In order to remedy this, we appeal to the polyfold technology of Hofer-Wysocki-Zehnder, which has been worked out in full detail in the Gromov-Witten case (i.e.\ for closed spheres) in \cite{HWZ_GW}.
In particular, before giving the proof, we need to recall the set up of the polyfold perturbation that we will use in order to prove our statement. The argument follows closely that of \cite[Section 7.2]{MNW}.

More precisely, we view $\overline{\mathcal{M}}_*$ as lying inside a \emph{Gromov-Witten polyfold} $\mathcal{B}_*$ \cite[Section\ 2.2, Definition\ 2.29, Section\ 3.5]{HWZ_GW} consisting of (not necessarily holomorphic) stable nodal configurations of spheres with one marked point and possibly multiple components. 
It comes with a natural evaluation map $ev:\mathcal{B}_*\rightarrow W_{cap}$, which extends the one of $\overline{\mathcal{M}}_*$. 
We view the nonlinear Cauchy-Riemann operator $\overline{\partial}_J$ as a \emph{Fredholm section} \cite[Definition\ 4.1]{HWZIII} of a \emph{strong polyfold bundle} $\mathcal{E}_*\rightarrow \mathcal{B}_*$ \cite[Definition\ 2.37, 2.38, Section 3.6]{HWZ_GW} with zero set $\overline{\partial}_J^{-1}(0)=\overline{\mathcal{M}}_*$.
We also have a forgetful map $\pi: \mathcal{B}_* \rightarrow \mathcal{B}$, where similarly $\mathcal{B}$ is a polyfold containing $\overline{\mathcal{M}}$, given by forgetting the marked point, and a strong polyfold bundle $\mathcal{E}\rightarrow \mathcal{B}$ satisfying $\pi^*\mathcal{E}=\mathcal{E}_*$, with a corresponding Fredholm section $\overline{\partial}_J$ making the obvious diagram commute.  

Since $\overline{\mathcal{M}}^\partial:=\{u_{t,q}\}\subset \overline{\mathcal{M}}$ consisting of the spheres in $(-\epsilon,0]\times \partial X \times \mathbb{S}^2$ is transversely cut out, the section $\overline{\partial}_J$ is in \emph{good position} \cite[Definition\ 4.12]{HWZIII} at the boundary. According to \cite[Theorem\ 4.22]{HWZIII}, we may introduce an abstract perturbation $\mathbf{p}$, which is a multivalued section of $\mathcal{E}_*$ \cite[Definition\ 3.35, Definition\ 3.43]{HWZIII}, so that:
\begin{itemize}
    \item[$\bullet$] $\overline{\partial}_J+\mathbf{p}$ is transverse to the zero section of $\mathcal{E}_*$;
    \item[$\bullet$] The perturbed moduli space $$\overline{\mathcal{M}}_*^\mathbf{p}=(\overline{\partial}_J+\mathbf{p})^{-1}(0)\subset \mathcal{B}_*$$ is a $(2n+2)$--dimensional compact, oriented, weighted branched orbifold with boundary and corners \cite[Section\ 3.2, Definition\ 3.22]{HWZIII}; and
    \item[$\bullet$] The perturbation $\mathbf{p}$ is supported away from a closed neighbourhood (in $\mathcal{B}_*$) of $\overline{\mathcal{M}}^\partial_*=\pi^{-1}(\overline{\mathcal{M}}^\partial)$.
\end{itemize}

\noindent We then define $\overline{\mathcal{M}}^\mathbf{p}\subset \mathcal{B}$ to be the image of $\overline{\mathcal{M}}^\mathbf{p}_*$ under the forgetful map $\pi^\mathbf{p}:=\pi\vert_{\calMbar^\mathbf{p}_*}$, which is a perturbation of $\overline{\mathcal{M}}$. 
Observe that the last condition implies that a collar neighbourhood of the boundary of $\overline{\mathcal{M}}^\mathbf{p}$ still consists of the holomorphic spheres from $\overline{\mathcal{M}}^\partial$, and local uniqueness still holds. 

Moreover, for $\mathbf{p}$ sufficiently small, every element in $\calMbar^\mathbf{p}$ will be close to {\em some} element in $\calMbar$, in the Gromov-Hausdorff topology. Hence we can assume that all elements in $\calMbar^\mathbf{p}$ are transverse to the co-cores and have precisely one intersection point, since this intersection is purely topological because of codimension reasons. 
Thus we obtain a well defined and smooth map $\mathcal{I}^\mathbf{p}: \calMbar^\mathbf{p}\to C'_+$, that agrees with $\mathcal{I}$ on $\overline{\mathcal{M}}^\partial$. 
\\

We are finally in place to prove \Cref{thm:homology_injection_boundary} in the non semi-positive case:
\begin{proof}[Proof of \Cref{thm:homology_injection_boundary} in the general case]
Fix any sufficiently small abstract perturbation $\mathbf{p}$ as described above, and consider the map $\mathcal{I}^\mathbf{p}_\pm\co \calMbar^\mathbf{p}\to C_\pm$. %\label{rem:fundamental_class}
Recall also that in this setting one has a well defined notion of integrating over the moduli space: we refer to \cite{HWZ_int,HWZ_GW} for further details. For our purposes it suffices to use the fact that there is a notion of \emph{sc-smooth differential forms} \cite[Definition 1.8]{HWZ_int} in $\overline{\mathcal{M}}^\mathbf{p}_*$, such that Stokes' theorem holds \cite[Theorem 1.11]{HWZ_int}, as well as a \emph{de--Rham cohomology group} $H^*_{dR}(\overline{\mathcal{M}}^\mathbf{p}_*)$ \cite[page 10]{HWZ_int}, such that said differential forms can be pulled back under the evaluation map. In our setting, one then has a well-defined notion of degree of a map, which agrees with the usual notion at the boundary of the moduli space, where smoothness is built into the construction. We then formally define (relative) homology $$H_*(\overline{\mathcal{M}}_*^\mathbf{p},\partial \overline{\mathcal{M}}^\mathbf{p}_*):=\left(H_{dR}^*(\overline{\mathcal{M}}^\mathbf{p}_*)\right)^*,$$ as the dual (over $\mathbb{R}$) of the de--Rham cohomology group, with the obvious notions of push-forwards, i.e.\ as duals of pull-backs.

This allows us to define a shriek map on homology (with real coefficients). Formally one does this by defining maps on the duals of the associated de Rham cohomologies via the integration map, as a map $ev^!: H_*(W_{cap}, \partial W_{cap} ) \to H_*(\overline{\mathcal{M}}_*^\mathbf{p},\partial \overline{\mathcal{M}}^\mathbf{p}_*)$ given via the equation:
$$ev^!(a)(\beta) = \int_{\calMbar^\mathbf{p}_{*}} ev^*PD(a) \wedge \beta, $$
where we have first used Poincar\'{e}-Lefschetz duality in the target. 

The remainder of the argument then goes through exactly as in the semi-positive setting and the theorem follows in the general case. 
More precisely, consider a rational cycle $\sigma$ in $N$ that is sent to $0$ via the map induced in homology by the inclusion $i\colon N\to \Wcap'$ as defined in  \Cref{eqn:diagram_pf_homol_inject}. 
Consider  $\calMstarbar^\mathbf{p}$, together with the perturbed intersection map $(\calI_\pm^\mathbf{p})_*$. 
Then, if $i_*[\sigma]$ is boundary of a homology class $b$ in $\Wcap$,
one can take its shriek $ev^!(b)$, seen as an element in the dual $(H^*_{dR}(\calMstarbar^\textbf{p}))^*$, and push it, via the map induced in the dual of the cohomology by $\calI_*^{\textbf{p}}$, to a class $c$ in the dual $(H^*_{dR}(C_\pm))^*$.
By Poincaré duality on $C_\pm$, we now have that $c$ can be seen as a class in $H_*(C_\pm;\QQ)$.
Moreover, as $ev$ at the boundary is actually a smooth map and, moreover, a diffeomorphism, we have $\partial c = (f_\pm)_*(\sigma)$.
But this means that $(f_\pm)_*[\sigma]=0$ in $H_*(C_\pm;\QQ)$, thus concluding the proof as in the semi-positive case.
%contradicts the fact that $f$ is homotopic to the standard inclusion of $\Sigma$ in $C'_+=\Sigma\times D^*\SSS^1$ as the first factor.
\end{proof}

\begin{remark} As an alternative to the above argument, one could appeal to the following. According to \cite[Remark 15.8]{HWZ_poly_fred} (cf.\ also \cite{McDuff_Fund_Class}) any compact oriented weighted branched orbifold admits a rational relative fundamental class. 
This then suffices to replicate the argument from the semi-positive case just using continuity of the maps, without the need to talk about sc-smoothness. Note, however, that our definition of relative homology bypasses the existence of this relative fundamental class.
\end{remark}

\medskip\noindent
\subsection{\textbf{The case of Bourgeois structures}}
\label{sec:bourgeois_case}

Let $BO(\Sigma,\phi)$ be the Bourgeois manifold obtained from the abstract contact open book $(\Sigma,\phi)$.
Denote more precisely $\beta=\alpha + \Phi_1 dx - \Phi_2 dy$ the contact form on $M\times \TT^2$ as in \Cref{sec:bourgeois_constr}, where $(x,y)$ are coordinates on $\TT^2$, $\Phi=(\Phi_1,\Phi_2)\colon M \to \RR^2$ is a (properly chosen) function describing the open book decomposition $M=OBD(\Sigma,\phi)$ and $\alpha$ is a contact form on $M$ supported by it.

Now, the vector field $\partial_y$ corresponding to the second factor of $\TT^2=\SSS^1_x\times\SSS^1_y$ of $M\times\TT^2$ is contact.
In particular, as described more generally in \Cref{sec_S_inv} for $\SSS^1-$equivariant contact structures, each hypersurface $M\times\SSS^1_x\times\{y_0\}$ is convex, and has hence a natural splitting as $V_+\cup \overline{V}_{-}$ as union of ideal Liouville domains.
In fact, $V_\pm$ are smoothly just given by the product of half of the open book $OBD(\Sigma,\phi)$ with $\SSS^1_x$, i.e.\ $\Sigma\times D^*\SSS^1_x$ with corners smoothened. 
The Bourgeois contact manifold is then nothing else than the gluing of the two Giroux domains $V_+\times\SSS^1_y$ and $V_-\times\SSS^1_y$ along their boundary.
(This decomposition is also explicitly described in \cite[Section 5.3]{DinGei12}, without the use of this terminology.)

\begin{remark}
\label{rmk:dual_viewpoint}
%As already said above, 
Since Bourgeois contact structures are $\TT^2$-invariant and we are free to choose any coordinate system $(x,y)$ on the $\TT^2$ factor, we can assume that $\SSS^1_x$ represents any primitive homology class in $H_1(\TT^2,\mathbb{Z})$.
%Here we considered $BO(\Sigma,\phi)$ as an $\SSS^1_y-$invariant contact structure, but one can also regard it as $\SSS^1_x-$invariant, and every construction that follows can be done using this other point of view.
This observation will be needed below in proving part of \Cref{thm:homology_injection_page}, namely that $\TT^2$ injects in integer homology in the filling $W$.
\end{remark}

We are now ready to prove \Cref{thm:homology_injection_page}.

\begin{proof}[Proof of \Cref{thm:homology_injection_page}]
 Since the inclusion of $\Sigma \subseteq \partial V_\pm=N$ into $V_{\pm} \cong \Sigma \times D^*\SSS^1$ is homotopic to the inclusion of a fiber, it is injective on homology and the injectivity of the map induced by $\Sigma\hookrightarrow W$ on rational homology follows directly from \Cref{thm:homology_injection_boundary}. 
 Likewise, in the case of semi-positive filling,  surjectivity at the level of fundamental groups of the inclusion $\Sigma\times\TT^2\hookrightarrow W$ also follows directly from \Cref{thm:homology_injection_boundary}.
 
 Furthermore, \Cref{thm:homology_injection_boundary} in fact also implies that the inclusion of the $\SSS^1_x$-factor is injective in rational homology. 
 Now taking $\SSS^1_x$ to represent an arbitrary (primitive) class in homology (c.f.\ \Cref{rmk:dual_viewpoint}), we can then conclude that the whole first homology group $H_1(\TT^2;\QQ)$ injects in $H_1(W;\QQ)$. 
 Furthermore, the statement that $H_*(\TT^2;\QQ)$ injects into $H_*(W;\QQ)$ is equivalent to the dual statement that the inclusion induces a surjection on cohomology $H^*(W,\QQ) \to H^*(\TT^2;\QQ)$ by the Universal Coefficient Theorem. 
 Using the cup-product structure on $H^*(\TT^2;\QQ)$, this follows from the injectivity in $H_1$,  concluding the proof.
\end{proof}

\subsection{\textbf{Applications to fillability}}\label{sec:applications}
We now prove those corollaries of \Cref{thm:homology_injection_page} that have been stated without proof in the Introduction. We start with \Cref{cor:stein_strong_equiv_homol_spheres} on the equivalence of strong symplectic and Stein fillability for Bourgeois manifolds associated to rational homology $3-$spheres:
\begin{proof}[Proof of \Cref{cor:stein_strong_equiv_homol_spheres}]
If $V=OBD(\Sigma,\phi)$ is a rational homology $3-$sphere, any inclusion of $\Sigma$ in the $3-$manifold $V$ as a page of the given open book is null homologous with coefficients $\QQ$.
In particular, \Cref{thm:homology_injection_page} tells that $BO(\Sigma,\phi)$ cannot be strongly fillable whenever $H_1(\Sigma;\QQ)$ is non-trivial, i.e.\ whenever $\Sigma$ is not a disc.
On the other hand, we know by \cite[Theorem A.(b)]{LMN} that $BO(\DD^2,\Id)$ is Stein fillable.
\end{proof}

We now prove \Cref{cor:pos_dehn_twist_no_asph_fill}, which states that Bourgeois $5-$manifolds associated to $3-$dimensional open books with planar pages and monodromies given by products of Dehn twists all of the same sign are weakly but not strongly symplectically fillable.
\begin{proof}[Proof of \Cref{cor:pos_dehn_twist_no_asph_fill}]
By \cite[Theorems A.(a) and B]{LMN} these contact structures are weakly fillable and, hence, it suffices to prove that these are not strongly fillable.
%According to \cite[Theorem B]{LMN}, it's enough to prove the case of the product of positive Dehn twists.
%That the Bourgeois manifold is weakly fillable then follows from the fact that the initial contact manifold is Stein fillable.
%
According to \Cref{thm:homology_injection_page}, it is hence enough to prove that the first Betti number of the manifold $b_1(OBD(\Sigma,\phi))$ is strictly smaller than that of the page  $b_1(\Sigma) = \# \partial \Sigma - 1$.  This can be seen as follows. 
Since for $\Sigma=\mathbb{D}^2$, the monodromy is necessarily trivial, we may assume that $\# \partial \Sigma\geq 2$. 
We can then attach $b_1(\Sigma)-1\geq 0$ discs to $\Sigma$, in order to cap all but $2$ of its boundary components; 
as $\phi$ is a product of Dehn twist of all the same signs, this can moreover be done in such a way that the resulting monodromy $\widehat{\phi}$ in the resulting annular page $\widehat{\Sigma} = D^*\SSS^1$ is a power $\tau^{\pm K}$ of the positive Dehn twist $\tau$, with $K>0$. 

Moreover, this $2-$dimensional $2-$handle attachment on the page can be realized by attaching a $4-$dimensional $2-$handle on the boundary component $OBD(\Sigma,\phi)\times\{1\}$ of a thickening $OBD(\Sigma,\phi)\times[0,1]$, along neighborhoods, of the form $\SSS^1 \times \DD^2$, of all but $2$ binding components. 
The resulting positive boundary of these $b_1(\Sigma)-1$ handle attachments is just $OBD(D^*\SSS^1,\tau^{\pm K})$ where $\tau$ is a Dehn twist and $K > 0$ as above. 
This manifold is then a Lens space and hence has $b_1=0$. 
Since attaching a $4-$dimensional $2$-handle reduces the first Betti number at the boundary by at most one, we deduce that $ b_1(OBD(\Sigma,\phi))\leq  b_1(\Sigma)-1$, as desired.
%Alternatively, one can attache handles of the form $\DD^2 \times D^*\TT^2$ along the product of the binding components with $\TT^2$ to obtain a strong symplectic cobordism to $BO(D^*\SSS^1,\tau^{\pm K})$, which is not strongly fillable again by \Cref{thm:homology_injection_page}. \textcolor{red}{I think this last point of view is very hard to understand now that the spine removal thing has been removed; I vote to remove it.}
\end{proof}

Now, we give a proof of \Cref{cor:commutator}, i.e.\ that if $\Sigma$ is a planar surface and $\BO(\Sigma,\phi)$ is strongly fillable, then $\phi$ lies in the commutator subgroup of the relative mapping class group:

\begin{proof}[Proof of \Cref{cor:commutator}]
We identify the pure braid group with the mapping class group of the punctured surface obtained by collapsing all but one boundary component of $\Sigma$ to a puncture. It is well known that  abelianization of the pure braid group on $n$-strands is the free abelian group generated by $\binom{n-1}{2}$ Dehn twists around pairs of punctures, where $n = \# \partial \Sigma$ is the number of boundary components (see for instance \cite[Chapter 9 and page 264]{FM}). 

Also, the kernel of the map on mapping class groups induced by the collapsing map described above is generated by the $n$ boundary parallel Dehn twists (see \cite[Theorem 3.18 and Proposition 3.19]{FM}). 
Hence the abelianization of $MCG(\Sigma,\partial \Sigma)$ is generated by $n$ boundary parallel Dehn twists and $\binom{n-1}{2}$ Dehn twists about pairs of boundary components different from a distinguished one; in total, it is hence generated by $\binom{n}{2} = n + \binom{n-1}{2}$ Dehn twists and has rank at most $\binom{n}{2}$.
%By collapsing out different boundary components, one deduces that this abelianization has rank $\binom{n}{2} = n + \binom{n-1}{2}$.
%One deduces that this abelianization is generated by $\binom{n}{2} = n + \binom{n-1}{2}$ Dehn twists.

Denote by $\widehat{\Sigma}_{i,j}$ the annulus obtained by capping off all but $2$ boundary components $\gamma_i,\gamma_j$ of $\Sigma$. 
Doing so for each pair of $i,j$ with $i\neq j$, one obtains a homomorphism
\begin{equation}
\label{eqn:MCG}
    MCG(\Sigma,\partial \Sigma) \to \prod_{i\neq j} MCG(\widehat{\Sigma}_{i,j},\gamma_i \cup \gamma_j) = \mathbb{Z}^{\binom{n}{2}},
\end{equation}
given by extending the element in the mapping class group of $\Sigma$ relative $\partial \Sigma$ as $\Id$ over the caps.
It is easily seen that this homomorphism is surjective.
Moreover, by the explicit description of the $\binom{n}{2}$ generators for the abelianization of $MCG(\Sigma,\partial \Sigma)$ above, it induces an injective homomorphism after taking the abelianization.
In particular, the map in \Cref{eqn:MCG} induces an isomorphism after abelianizing.

Hence, if the monodromy is not in the commutator subgroup, it is mapped to a non-trivial element in $\mathbb{Z}^{\binom{n}{2}}$, and so it is non-trivial in $MCG(\widehat{\Sigma}_{i_{0},j_{0}},\gamma_{i_0} \cup \gamma_{j_0})\cong \mathbb{Z}$ after capping off all but $2$ well-chosen boundary components $\gamma_{i_0}, \gamma_{j_0}$ of the page. 
We then deduce as in the proof of \Cref{cor:pos_dehn_twist_no_asph_fill} that $b_1(OBD(\Sigma,\phi)) < b_1(\Sigma)$, and we conclude the argument by appealing to \Cref{thm:homology_injection_page}.
\end{proof}

We now prove \Cref{cor:high_dim_stab_not_fill} about the fillability of Bourgeois structures arising from open books which are positive stabilizations. 

\begin{proof}[Proof of \Cref{cor:high_dim_stab_not_fill}]
The abstract positive stabilization $OBD(\Sigma_+,\phi_+)$ of a contact open book $OBD(\Sigma,\phi)=(M,\xi)$ can be interpreted as the open book on $M\#\SSS^{2n+1}$ given by the Murasugi sum of $OBD(\Sigma,\phi)$ on $M$ with the open book $OBD(D^*\SSS^n,\tau)=(\SSS^{2n+1},\xi_{std})$, where $\tau$ is the positive Dehn-Seidel twist.
More precisely, from the explicit construction of Murasugi sum (see for instance \cite[Proposition 2.6]{CasMur16}), one can see that %, I CHANGED THIS on 25 Mai at 18:00 for a properly chosen page $\Sigma$, 
the additional Lagrangian sphere $S$ created in the stabilized page $\Sigma_+$ is actually contained in the contact ball $\DD^{2n+1}=\SSS^{2n+1}\backslash D$, where $D$ is the small Darboux ball taken out of the $\SSS^{2n+1}$ factor to build the connected sum $M\#\,\SSS^{2n+1}$. 
Hence, $S$ is null-homologous in $M\#\,\SSS^{2n+1}=M$.
In particular, as $S$ is not null-homologous (by construction) in $\Sigma_+$, \Cref{thm:homology_injection_page} gives that $BO(\Sigma_+,\phi_+)$ is not strongly fillable.
%Consider the co-core $C$ of the handle $H$ such that $\Sigma_+=\Sigma \cup H$, seen as lying in the page at angle $\theta=0\in \RR/\ZZ$ in the open book $OBD(\Sigma_+,\phi_+)$.
%Consider now any connection on $OBD(\Sigma_+,\phi_+)\setminus B \to \SSS^1$ whose time$-1$ flow sends pages to pages and gives the monodromy $\phi_+$; denote $\psi_t$ its flow.
%The union of the images $\psi_t(C)$ for $t\in[0,1$ then gives a $(n+1)-$chain in the (integer) singular chain complex of $OBD(\Sigma_+,\phi_+)$, whose boundary is just the $n-$cycle $C + \tau(C)$, where $\tau$ is the positive Dehn twist along the Lagrangian sphere $S$ coming from the positive stabilization.
%Now, $C+\tau(C)$ is just cohomologous to $\pm [S]$ (here, the sign depends on the action of $\tau$ on the fundamental class of the Lagrangian sphere, i.e.\ on the parity of $n$).
%In particular, $[S]=0\in H_*(OBD(\Sigma_+,\phi_+);\ZZ)$.
%As $[S]\neq0\in H_*(\Sigma;\ZZ)$, \Cref{thm:homology_injection_page} then implies that $BO(\Sigma_+,\phi_+)$ is not strongly fillable, as desired.
%\textcolor{red}{Should work, please double check}
\end{proof}

\noindent We now give a proof of  \Cref{thm:BOfill} on the fillability of $BO(D^*\SSS^n,\tau^k)$, of which \Cref{thm:bou_tau_high_dim} in the case $k=1$ is an immediate consequence:

\begin{proof}[Proof of \Cref{thm:BOfill}]
It follows from \Cref{thm:strong_cobord_bourgeois} and \cite[Theorems A.(b) and B]{LMN} that $BOFill(n)$ is a subgroup of $\mathbb{Z}$.
It is hence necessarily cyclic.
%Note that $BO(D^*\SSS^n,\mathbb{\tau}^k)\cong BO(D^*\SSS^n,\mathbb{\tau}^{-k})$ as contact manifolds from \cite[Theorem B]{LMN}, and so it suffices to assume $k\geq 0$.

According to \Cref{thm:homology_injection_page}, in order to determine the generator of $BOFill(n)\le \ZZ$, it then suffices to check for which $k\geq 1$ the homology of the page injects in $\Sigma_{k,n}:=OBD(D^*\SSS^n,\tau^k)$.
%(Recall that $BO(\Sigma,\Id)$ is Stein fillable \cite[Theorem A.(b)]{LMN}.)
%In the case $n=1$, we have $\Sigma_{k,1}=L(k,k-1)$ if $k\geq 1$, and
%$\Sigma_{0,1}=\SSS^2\times\SSS^1$.
%%$\Sigma_{k,1}=L(k,1)$ if $k\leq 0$. 
%If $k\neq 0$, these Lens spaces have vanishing rational first homology group. 
%We then conclude that $k_0(1)=0$. 
%More generally, 
According to \cite[Proposition 4.10]{KvK16}, for $n,k\geq 1$ the manifold $\Sigma_{k,n}$ is actually the $(2n+1)-$dimensional Brieskorn manifold $\Sigma(2,\dots,2,k)$. 
An explicit computation via Mayer-Vietoris, using the Heegaard splitting associated to the open book, gives the following: 
if $n$ is odd, then $H_n(\Sigma_{k,n};\mathbb{Z})=\mathbb{Z}_k$ for every $k$ and so it vanishes over $\mathbb{Q}$; 
and $H_n(\Sigma_{k,n};\mathbb{Z})=0$ if $n$ is even, and $k\geq 1$ is odd. 
We then see that $k_0(n)$ is even. 
Our method is inconclusive for $n$ and $k\geq 2$ both even, since we have $H_n(\Sigma_{k,n};\mathbb{Q})=\mathbb{Q}$, generated by the zero section of the page. Indeed, if $k$ is such that $\tau^k$ is smoothly isotopic to $\Id$, then $\Sigma_{k,n}\cong  \mathbb{S}^n\times\mathbb{S}^{n+1}= \Sigma_{0,n}$ as smooth manifolds, hence the conclusion of \Cref{thm:homology_injection_page} clearly holds, and we obtain no obstruction (according to \cite{Keating2021OnTO}, cf.\ \cite{K07,Cer70}, for arbitrary even $n\geq 2$, $\tau$ is known to have finite order as a smooth map). Alternatively, for odd $n$, one may observe that $\tau^k$ acts non-trivially in homology if $k\neq 0$ (by the Picard--Lefschetz formula), thus \Cref{rk:phi_is_id_homology} allows to conclude; however, $\tau^2$ is homologically trivial if $n$ is even, so this argument is inconclusive as well.
\end{proof}

%We thus obtain infinitely many non-diffeomorphic examples of weakly but not strongly fillable $2n+3$-manifolds, $n\geq 1$.

\section{\textbf{Symplectically aspherical fillings of Bourgeois Contact Manifolds}}\label{sec:fillings_Bour}

Let $W$ be a symplectically aspherical strong symplectic filling of $BO(\Sigma,\phi)$. The goal of this section is to prove \Cref{thm:hom_apsherical_filling_Bour} and \Cref{thm:application} from the Introduction. 

The symplectic manifold $W_{cap}$ resulting from the handle attachment described in \Cref{lemma:capping_cobord_high_dim} %(i.e.\ the one in \cite[Theorem 6.1]{MNW})
then has boundary of the form $N \times \SSS^2$, where $N$ is the contact boundary of $\Sigma \times D^*\SSS^1_x$ (smoothened at the corners).
More precisely, in a neighborhood $(-\delta,0]\times N \times \SSS^2$ of the boundary, the symplectic structure has normal form $d(e^t \alpha) + \omega_S$, with $\alpha$ a contact form on $N$ and $\omega_S$ an area form on $\SSS^2$. 

We thus obtain a moduli space $\mathcal{M}$ consisting of holomorphic spheres inside $W_{cap}$, arising from the $\SSS^2-$factor at the boundary. 
Under the asphericity assumption it turns out that there are no nodal degenerations in the Gromov compactification and we thus obtain:

\begin{proposition}\label{prop:compactnessTn}
    $\mathcal{M}$ is a compact oriented smooth manifold (with boundary). 
\end{proposition}
\begin{proof}

As the filling $W$ is symplectically aspherical, there cannot be any nodal curve in $\calMbar$ which has one closed sphere component not intersecting either of the co-cores.
Hence, the only nodal degeneration that could a priori occur is one separating a regular sphere in the moduli space in multiple spheres, each intersecting at least one co-core.
However, this would imply, after puncturing the spheres at their intersection points with the co-cores, that one of the $\SSS^1-$factors in the $\TT^2-$factor of $\partial W = OBD(\Sigma,\phi) \times \TT^2$ is null-homologous in $W$. 
(Here, we use the fact that removing the co-cores from $W_{cap}$ results in a manifold which deformation retracts onto $W$.)
This would then contradict \Cref{thm:homology_injection_page}, so that no nodal degenerations are possible, as desired.
\end{proof}

\noindent \textbf{A moduli space of punctured curves.} 
As in \Cref{sec:general_obstr_fillab}, we can consider the marked moduli space $\mathcal{M}_*$ obtained by adding a marked point to the domains, which comes equipped with an evaluation map $ev:\mathcal{M}_*\rightarrow W_{cap}$, and a forgetful map $\pi:\mathcal{M}_*\rightarrow\mathcal{M}$. 
%Similarly, we have the same structure for its Gromov compactification $\overline{\mathcal{M}}$, whose elements intersect the $J$-invariant co-cores precisely once. 

We now consider a punctured version of the compact marked moduli space $\calM_*$.
More precisely, as each sphere in $\calM$ intersects the co-cores $C_\pm$ positively and transversely in exactly one point, we may remove from each sphere in $\calM_*$ small disc-like neighborhoods around its intersection with the two co-cores of \Cref{lemma:capping_cobord_high_dim}.
This results in a (compact) marked moduli space $\calM_*^{c}$ of (compact) cylinders.
Moreover, this has a natural fibration structure $\pi$ over $\calM^{c}\coloneqq\calM$ induced by the forgetful map $\calM_*\to\calM$, and an evaluation map $ev\colon \calM_*^{c}\to W$, where we see here $\Wcap$ minus a small tubular neighborhood of the co-cores as naturally diffeomorphic to $W$. 
Because of the foliation property of the moduli space of spheres $\calM_*$ near $\partial \Wcap$, the evaluation map of the cylindrical counterpart $\calM_*^{c}$ is a diffeomorphism from a neighborhood of the vertical boundary $\partial_v\calM_*^{c} \coloneqq \pi^{-1}(\partial \calM)$ to a neighborhood of 
$N \times D^*\SSS^1\subset BO(\Sigma,\phi)=\partial W$ inside $W$, where $N$ is the boundary of $V_\pm=\Sigma\times D^*\SSS^1_x$ with corners smoothened as in the description in the beginning of \Cref{sec:bourgeois_case}.
%$\Sigma \times D^*\SSS^1\subset BO(\Sigma,\phi)=\partial W$ inside $W$.
Hence, the evaluation map has in particular degree $1$.

What is more, it also sends the horizontal boundary $\partial_h\calM_*^{c} \coloneqq \overline{ \partial\calM_*^{c} \setminus \partial_v\calM_*^{c}}$ into a small neighbourhood of the co-cores minus the co-cores themselves, i.e.\ into $N_\delta(C_\pm)\setminus C_\pm$.
In particular, each cylinder in $\calM^{c}_*$ has boundary components which have a natural sign $\pm$ according to where they are mapped to via $ev$; this gives a partition $\partial_h\calM^{c}_* = \partial_h^+ \calM^{c}_* \sqcup \partial_h^- \calM^{c}_* $. 
We further note that the image of each of these components retracts onto the piece $V_\pm \times \SSS^1$ of $BO(\Sigma,\phi)=\partial W$ in the decomposition induced by the $\SSS^1$-invariant perspective (cf.\ \Cref{sec_S_inv}).

\medskip

\noindent \textbf{Homology of aspherical fillings.}
We now prove that the inclusion of each $V_{\pm} \times \SSS^1 = \Sigma \times D^*\SSS^1 \times \SSS^1 \simeq \Sigma \times \SSS^1 \times \SSS^1$ into an aspherical filling induces an isomorphism in integral homology. 

\begin{proof}[Proof of \Cref{thm:hom_apsherical_filling_Bour}]
	Recall that Lefschetz duality states that the cap product with the fundamental class living in the top degree homology relative to the boundary induces an isomorphism between relative cohomology groups and homology groups (with proper degree shifting). 
	Hence, since the evaluation map $ev$ of $\calM_*^{c}$ has degree $1$, naturality of the cap product implies that it induces a surjection on homology with $\ZZ$-coefficients. 
	Next, note that the inclusion  $\partial_h^\pm \calM^{c}_* \to \calM^{c}_*$ induces a homotopy equivalence, as each cylinder can be smoothly collapsed to either of its boundary components.
	Now, the image of $\partial_h^\pm \calM^{c}_*$ under $ev$ can be homotoped into $V_\pm \times \SSS^1$. 
	Then, 
	one can do the following: 
	given any homology class in $W$, one can pull it back 
	%(after perturbation) 
	to the moduli space, homotope it
	%smoothly 
	to the positive boundary $\partial_h^+\mathcal{M}^c$, and finally map it near $V_\pm\times\SSS^1$ 
	%(up to homotopy) 
	via $ev$. 
	We hence obtain surjectivity of the inclusion $V_\pm\times\SSS^1\hookrightarrow W$ as stated in \Cref{thm:hom_apsherical_filling_Bour}.
	
	For injectivity, we note that the inclusion of  $\Sigma \times \partial_+(D^*\SSS^1_x) \times \SSS^1_y \subseteq \partial V_+\times \SSS^1_y$ into $V_+ \times \SSS^1_y$ is a homotopy equivalence, where $\partial_+(D^*\SSS^1)$ denotes the positive boundary component of $D^*\SSS^1$. Hence any non-trivial homology class $[c] \ne 0$ in $V_+ \times \SSS^1$ is represented by a class in $\partial V_+ \times \SSS^1$, along which the evaluation map is a diffeomorphism. Thus if $c = \partial b$ were the boundary of a singular chain, we could pull back, after perturbing the evaluation map (relative to the horizontal boundary), and then push the resulting chain $\bar{b} = ev^{-1}b$ into $\partial_h^+ \calM^{c}_*$. Pushing forward by $ev$ we would obtain that $c$ is a boundary of $ev_*\bar{b}$ in $ V_+ \times \SSS^1$, which is a contradiction.
	The case of $V_-$ is completely analogous.
	
	Finally observe that the inclusion of $\Sigma \times \partial_+(D^*\SSS^1_x) \times \SSS^1_y$ into $\partial W$ is precisely the inclusion of $\Sigma \times \TT^2 \hookrightarrow OBD(\Sigma,\phi) \times \TT^2$, concluding the proof.
\end{proof}
%One can now factor (up to homotopy) $j\co \TT^{n} \hookrightarrow W$ from the statement of \Cref{prop:iso_homol_fundgp_fill_STn} as composition of a map $i\co \TT^n \to \mathcal{M}_*^{c}$ and the evaluation map $ev\co \mathcal{M}_*^{c} \to W$, where $i$ is defined as follows. First, note that the $j_0$ defined above can be homotoped to the inclusion $j_0'$ of $\TT^n$ as $\{p_0\}\times \TT^{n-1}\times \{+1\}\times \SSS^1$ in the positive component $D^*\TT^{n-1}\times\{+1\}\times \SSS^1$ of $D^*\TT^{n-1}\times S^*\SSS^1$, where $p_0$ is a point in the cotangent direction of $D^*\TT^{n-1}$ in the boundary of the cotangent fiber.
%In particular, this inclusion actually lives in the intersection of the two components $D^*\TT^{n-1}\times S^*\SSS^1$ and $S^*\TT^{n-1}\times D^*\SSS^1$; here, $ev$ is a diffeomorphism \textcolor{red}{We should remove this reference to the paper I think}(as it is near the paper) and this inclusion can then be pulled back to $\calM_*^{c}$.
%The composition of $j_0'$ with the inverse of $ev$ then gives the desired map $i$.

%The fact that the moduli space is smooth then allows for the following strengthening  of \Cref{thm:homology_injection_boundary} to integral homology, since in this case one can avoid pseudo-cycles.

%ADD

\subsection{\textbf{Symplectically aspherical fillings of $S^*\mathbb{T}^n$}}\label{sec:fillings_STn}

In this section, we consider $(S^*\mathbb{T}^n,\xi_{std})$, the unit cotangent bundle of $\mathbb{T}^n,$ for $n\geq 2$, with its standard Stein fillable contact structure, and we prove \Cref{thm:application} from the Introduction. 

The Bourgeois contact manifold given by $ BO(D^*\TT^{n-2},\Id)$ is contactomorphic to the contact boundary of $D^*\TT^{n-2} \times D^*\TT^{2}$ by \cite[Theorem A.(b)]{LMN}, which can in turn be identified with the boundary of of the unit cotangent bundle of the $n$-torus namely $(S^*\TT^n, \xi_{st})$. 
In particular, $\xi_{std}$ is an $\SSS^1-$invariant contact structure, with respect to the second $\SSS^1$-factor of the Bourgeois torus $\TT^2$. %, as constructed in \cite[Section 2]{GeiSti10}. 

%\textbf{Properties of $\calMbar$.}
%We first notice that the moduli space $\calMbar$ gives a foliation by spheres near the boundary by \Cref{lemma:localuniq}. We now need to understand degenerations in $\calMbar$. To this end we start with the following topological statement, which follows from the proof of \Cref{thm:homology_injection_page}:
%, using the moduli space $\overline{\mathcal{M}}$ coming from two different spinal open books as in \Cref{eqn:spinal_obd_SstarTn}, say $j=1$ and $j=n-1$, one gets:
%\begin{proposition}\label{prop:Tnsurvives}
 %  The natural embedding of a section $\mathbb{T}^{n}$ of $\partial W = S^*\mathbb{T}^{n} \to \TT^{n}$ into $W$ is $H_*-$injective.
   %$\subset \partial W_{cap}$ survives homologically in $W_{cap}$.
%\end{proposition}

%\begin{proof}
%It follows directly from the proof of \Cref{thm:homology_injection_page} applied to the various $\SSS^1$-fibrations on $S^*\TT^n$, that the inclusion of any of the $\SSS^1$-factors survives in the filling. 
%Dually one then has degree $1$ cohomology classes $\gamma_i$ on the filling that evaluate non-trivially on the $i-$th circle factor $\SSS^1_i$. 
%Then using the cup-product structure on $\TT^n$, it follows that $\gamma_1 \smile  \dots \smile \gamma_n$ restricts to a non-trivial class in $H^n(\TT^n)$. 
%This implies that the induced map on cohomology is surjective and hence the induced map on homology is injective, as the homology groups of $\TT^n$ are torsion free.
%\end{proof}

\medskip 

\noindent \textbf{Symplectically aspherical fillings of $S^*\TT^n$: homotopy type.}
%%As a first step towards \Cref{thm:application}, we will prove that $W$ is homotopy equivalent to $D^*\TT^n$.
%% Let $W$ be a symplectically aspherical filling of $S^*\TT^n = BO(D^*\TT^{n-2},\Id)$.
%We next consider $S^*\TT^n$ as the contact boundary of the product $D^*\TT^{n-1}\times D^*\SSS^{1}$, which then decomposes as
%\begin{equation}
%\label{eqn:dec_SstarSn}
%S^*\TT^n=D^*\TT^{n-1}\times S^*\SSS^1 \bigcup S^*\TT^{n-1}\times D^*\SSS^1,
%\end{equation}
%where the first (disconnected) piece $D^*\TT^{n-1}\times S^*\SSS^1$ is the product of the union of the two Liouville manifolds (both $D^*\TT^{n-1}$, 
%but one with orientation reversed) with $\SSS^1$. %together with their natural contact forms, and the second piece $S^*\TT^{n-1}\times D^*\SSS^1$ is the region in which the interpolation between the contact forms on the other two pieces occurs.
%In particular, writing $S^*\SSS^1 = \{\pm1\}\times\SSS^1$, the piece $D^*\TT^{n-1}\times S^*\SSS^1$ has two connected components, namely $D^*\TT^{n-1}\times \{+1\}\times\SSS^1$, and $D^*\TT^{n-1}\times \{-1\}\times\SSS^1$, which we call the \emph{positive} and \emph{negative} components respectively.
%Notice that there is a natural map $j_0\co \TT^n \to S^*\TT^n$ which factors through the inclusion of $\TT^n=\{0\}\times\TT^{n-1}\times\{+1\}\times\SSS^1$ into the positive component of the component $D^*\TT^{n-1}\times S^*\SSS^1$ of $S^*\TT^n$.  
Recall that $S^*\TT^n = BO(D^*\TT^{n-2},\Id)$.
According to the $\SSS^1-$invariant picture as in \Cref{sec_S_inv}, this gives a decomposition $$S^*\TT^n = V_+\times\SSS^1\cup \overline{V}_-\times\SSS^1,$$ where $V_\pm$ are smoothly given by $D^*\TT^{n-2}\times D^*\SSS^1$, i.e.\  $ D^*\TT^{n-1}$ (up to rounding corners). 
Denote then by $j_0\colon \TT^n \hookrightarrow S^*\TT^n$ the inclusion given by the composition of inclusions $\TT^n \hookrightarrow D^*\TT^{n-1} \times \SSS^1\hookrightarrow S^*\TT^n$, where the first injection is just induced by the zero section $\TT^{n-1}\to D^*\TT^{n-1}$.

We next use the following result in order to conclude that the given aspherical filling $(W,\omega)$ is homotopy equivalent to $D^*\TT^n$. 
\begin{proposition}
	\label{prop:iso_homol_fundgp_fill_STn}
	Let $j \co \TT^{n}\hookrightarrow W$ be given by the composition of $j_0$ defined above and the natural inclusion $S^*\TT^n = \partial W$ into $W$. 
	Consider also a lift $\tildej$ of $j$ to the universal covers $\RR^n$ and $\tildeW$.
	Then:
	\begin{enumerate}
		\item $H_0(\tildeW;\ZZ)=\ZZ$ and $H_k(\tildeW;\ZZ)=\{0\}$ for $k>0$, 
		\item $j_*\co \pi_1(\TT^n)\to \pi_1(W)$ is an isomorphism.
	\end{enumerate}
\end{proposition}
According to \cite[Section 4.2, Exercise 12]{HatAlgTop}), this then implies that the inclusion $j$, and hence the inclusion of its (trivial) normal bundle $D^*\TT^n$, induces a homotopy equivalence, as desired. 

\begin{remark}
We point out that, as opposed to the simply connected case, it is not in general true that a map inducing an isomorphism on fundamental group and isomorphisms on homology is a homotopy equivalence. 
This is the reason why we need to consider the universal cover in \Cref{prop:iso_homol_fundgp_fill_STn}.
\end{remark}

\begin{proof}[Proof (\Cref{prop:iso_homol_fundgp_fill_STn})]
	We consider the moduli space of punctured curves and its evaluation map $\calM_*^{c} \to W$ considered as in the proof of \Cref{thm:hom_apsherical_filling_Bour}. Since the evaluation map $ev$ of $\calM_*^{c}$ has degree $1$, it induces a surjection on $\pi_1$. 
	Moreover, recall that the moduli space retracts onto its positive boundary, which is mapped onto a neighborhood of the positive component, that in this setting is just $D^*\TT^{n-1}\times \ \{+1\} \times \SSS^1$. 
	Arguing as in the proof of \Cref{thm:hom_apsherical_filling_Bour}, one concludes that $j\co \TT^n \hookrightarrow W$ also induces a surjection on fundamental group.
	In particular, we deduce that the fundamental group of the filling is abelian. 
	The fact that $j$ is $H_*-$injective follows from \Cref{thm:hom_apsherical_filling_Bour}; 
	%We now prove that $j$ is also $H_*$--injective;
	as $\pi_1(\TT^n)$ and $j_*(\pi_1(\TT^n))=\pi_1(W)$ are abelian, this immediately implies that $j$ is also $\pi_1$--injective.
%	Consider a class $x\in H_*(\TT^{n})$, and denote $y\coeq i_*(x)$ in $H_*(\calM_*^{c})$.
%	Suppose that $ev_*(y)$ is zero in $H_*(W)$;
%	this means that there is $z\in H_*(W)$ with $ev_*(y)=\partial z$. 
%	As $\calM_*^{c}$ and $ev$ are both smooth, there is a well defined $ev^!(z)$ in $H_*(\calM_*^{c})$ such that $ev_* ev^! (z) = z$, where $ev^!$ is given, geometrically, by perturbing $ev$ to be transverse to a cycle representing $z$ and taking its preimage.
%	
%	Moreover, as $y$ is supported near $\partial^+_h\calM_*^{c}$, where $ev$ restricts to a diffeomorphism, $\partial ev^!(z)=y$, so that $y=0$ in $H_*(\calM_*^{c})$. 
%	Moreover, as $\calM_*^{c}$ retracts onto $\partial^+_h\calM_*^{c}$, $y=0$ in $H_*(\partial^+_h\calM_*^{c})$ too.
%	Now, because $ev(\partial^+_h\calM_*^{c})=D^*\TT^{n-1}\times\SSS^1$ is the component $D^*\TT^{n-1}\times S^*\SSS^1$ of $S^*\TT^n=\partial W$, this means $j_*(x)=ev_*(y)=0$ in $H_*(D^*\TT^{n-1}\times\SSS^1)$, i.e.\ $x=0$ in $H_*(\TT^n)$, as $j$ factorizes as the natural inclusion $\TT^n=\TT^{n-1}\times\SSS^1\hookrightarrow D^*\TT^{n-1}\times\SSS^1 \hookrightarrow \partial W\hookrightarrow W$.
	
	The only thing left to prove is hence that the map $\tildej$ on the universal cover also induces an isomorphism in $H_*$. 
	Injectivity follows trivially from the fact that the universal cover of $\mathbb{T}^n$ is contractible.
	We then prove surjectivity.
	
	In order to do so, we first make the following observations. 
	All holomorphic curves in $\mathcal{M}^{c}$ naturally lift to maps from the plane to the universal cover $\widetilde{W}$ of $W$. 
    This gives a corresponding smooth moduli space $\widetilde{\calM}$ of infinite holomorphic strips in $\widetilde{W}$, together with its marked version $\widetilde{\calM}_*$ equipped with an evaluation map $\widetilde{ev}\co \widetilde{\calM}_{*} \to \widetilde{W}$.
    One can check that the $\widetilde{ev}$ is proper and that has degree $1$ (using cohomology with compact supports), since the diffeomorphism on the vertical boundary $\partial_v \calM^c_*$ of $\calM^c_*$ induced by $ev$ lifts to a diffeomorphism from a subset of the boundary of $\widetilde{\mathcal{M}}_*$ to a subset of the boundary of $\widetilde{W}$.
    %considered as a subset of the universal cover of $W$. 
    Notice now that $\widetilde{\calM}_{*}$ deformation retracts onto its positive horizontal boundary $\partial_{h}^+ \widetilde{\calM}_{*}$, defined as the lift of $\partial_h^+\calM_*^c$, i.e.\ as the union of all the lifts of the positive boundaries of the cylinders in $\calM_*^{c}$.
    As $j$ is an isomorphism in $\pi_1$, using standard covering space arguments, one can check that this set is mapped, via the lifted evaluation map, to a lift $D^*\RR^{n-1}\times \RR$ of the piece $D^*\TT^{n-1}\times\SSS^1$ of $\partial W \subset W$ to $\widetilde{W}$. 
    %An argument analogous to the one of the previous paragraph then allows to deduce that any homology class in $\widetilde{W}$ comes from $D^*\RR^{n-1}\times \RR$ so that the relative homology $H_*(\widetilde{W}, \partial \widetilde {W})$ is trivial in all degrees.
    
    Using this setup, we can now prove that any homology class in $\widetilde{W}$ comes from $D^*\RR^{n-1}\times \RR$, which then concludes the proof of \Cref{prop:iso_homol_fundgp_fill_STn}. 
    %; this is enough to conclude the proof of the Proposition, as this implies that the relative homology $H_*(\widetilde{W}, \partial \widetilde {W})$ is trivial in all degrees.
    More precisely, consider a class $x\in H_k(\widetilde{W},\ZZ)$, for $k\geq1$. 
    (The case $k=0$ simply follows from the fact that $\widetilde{W}$ is connected.)
	As $\widetilde{\calM}_*$ is an orientable manifold (with boundary) and $\widetilde{ev}$ has degree $1$, there is a well defined $\widetilde{ev}^!(x)$ in $H_*(\widetilde{\calM}_*)$ such that $\widetilde{ev}_* \widetilde{ev}^! (x) = x$, where $\widetilde{ev}^!$ is given, geometrically, by perturbing $\widetilde{ev}$ to be transverse to a cycle representing $x$ and taking its preimage.
	Now, as $\widetilde{\calM}_*$ retracts onto $\partial^+_h\widetilde{\calM}_*$, $\widetilde{ev}^!(x)$ can be homotoped to $\partial^+_h\widetilde{\calM}_*$.
	In particular, $x=\widetilde{ev}_* \widetilde{ev}^! (x)$ is homologous in $\widetilde{W}$ to a cycle $\sigma$ in $\widetilde{ev}(\partial^+_h\widetilde{\calM}_*)=D^*\RR^{n-1}\times\RR$.
	As the latter is contractible, $[\sigma]$, and so $x$, is null-homologous in $\widetilde{W}$, as desired.
\end{proof}

\noindent \textbf{Symplectically aspherical fillings of $S^*\TT^n$: diffeomorphism type.} 
Once the homotopy type is understood, the diffeomorphism type can be determined using the s--cobordism theorem. 
The argument below is just an adaptation of \cite[Sections 5 and 8]{BGZ} to our setting.
We thus give a sketch of the proof, referring to the proofs of the technical statements in \cite{BGZ}; for the readers' ease, we also adopt their notations.

We start by describing the spaces involved in the argument. 
Let $W_1$ be the result of attaching a topologically trivial cobordism $[0,1]\times S^*\TT^n$ to $W$ along its boundary $M_0\coeq S^*\TT^n = \{0\}\times S^*\TT^n$.
%Consider also on $S^*\TT^n$ the decomposition given by 
%$$
%S^*\TT^n=\partial D^*\mathbb{T}^n=D^*\TT^{n-1}\times S^*\SSS^1 \bigcup S^*\TT^{n-1}\times D^*\SSS^1.
%$$
%Recall that the piece $D^*\TT^{n-1}\times S^*\SSS^1$ has a ``positive'' component $D^*\TT^{n-1}\times \{+1\}\times\SSS^1$ given by the natural identification $S^*\SSS^{1}=\{\pm1\}\times \SSS^1\subset \RR\times\SSS^1$.
%The cobordism $[0,1]\times S^*\TT^n$ then contains a collar, diffeomorphic to $D^*\TT^{n-1}\times D^*\SSS^1=D^*\TT^n$, of the positive component of the piece $\{1/2\}\times D^*\TT^{n-1}\times S^*\SSS^1$ of $\{1/2\}\times S^*\TT^n$.

Recall from the previous subsection that we have a natural inclusion $j_0\colon \TT^n \to \partial W = S^*\TT^n$.
As $j_0(\TT^n)$ has trivial normal bundle in $\partial W$,
%In other words, 
there is a (smooth) copy of $W_0:=D^*\mathbb{T}^n$ which is entirely contained in the cobordism $[0,1]\times S^*\TT^n$.
Let then $X\coeq W_1\setminus W_0$ and $M_1\coeq \partial W_1$; notice that $\partial X = M_1 \cup (-M_0)$.
The aim is now to prove that $X$ is diffeomorphic to a cylinder $[0,1]\times S^*\TT^n$, so that $W_1$ is actually diffeomorphic to $W_0$, as desired.

As explained in \cite[Lemmas 5.1 and 5.2]{BGZ}, the fact that $W$ is homotopy equivalent to $D^*\TT^n$ implies that the inclusions $M_0,M_1\hookrightarrow X$ induce isomorphisms on $\pi_1$ and on $H_*$.
Moreover, as $S^*\TT^n$ is a \emph{simple space} (i.e.\ the action of its $\pi_1$ on every homotopy group is trivial), arguing exactly as in \cite[Lemmas 8.1 and 8.2]{BGZ} one can show that $M_0,M_1\hookrightarrow X$ actually induce isomorphisms on all homotopy groups.
This proves that $X$ is an h--cobordism between $M_0$ and $M_1$.

Now, as $M_0=S^*\TT^n$, the Whitehead group $\mathrm{Wh}(\pi_1(M_0))$ vanishes, so that the Whitehead torsion of the inclusion $M_0\hookrightarrow X$ is necessarily zero.
The s--cobordism theorem then tells that $X$ is diffeomorphic to $[0,1]\times S^*\TT^n$, as desired.

\section{\textbf{Further discussion and open questions}}  \label{sec:further_q}
Our results are, together with \cite{LMN}, among the first steps in understanding the nature of the contact structures given by Bourgeois' construction, and several open questions remain. Firstly:

\begin{question}
Are Bourgeois contact structures tight in \emph{all} odd dimensions?
Moreover, is every Bourgeois contact structure \emph{weakly} fillable, at least in dimension $5$?
\end{question} 

It is an important problem to understand more precisely the dependence of the Bourgeois structure on the starting open book decomposition. 
By a direct consequence of their definition, all Bourgeois contact structures are contact deformations of the almost contact structure $\xi_V\oplus T\mathbb{T}^2$ (i.e.\ the endpoint $\eta_1$ of a path $(\eta_t)_{t\in[0,1]}$ of hyperplane fields starting at $\eta_0=\xi_V\oplus T\mathbb{T}^2$ and such that $\eta_t$ is contact for $t>0$). 
One can then construct, as in \cite[Example 1.1]{MNW}, weak cobordisms between $BO(\Sigma,\phi)$ and $BO(\Sigma^\prime, \phi^\prime)$ for any $OBD(\Sigma^\prime, \phi^\prime)$ and $OBD(\Sigma, \phi)$ supporting the same contact structure. 
What's more, besides sharing the formal homotopy class, \Cref{thm:tight_bourgeois} in this paper shows in particular that the tight vs overtwisted classification type of any 5-dimensional $BO(\Sigma,\phi)$ is independent of the open book.

On the other hand, in \cite[Corollaries 10.6 and 10.8]{Bothesis}, Bourgeois used cylindrical contact homology, with respect to noncontractible homotopy classes of Reeb orbits, in order to distinguish infinitely many Bourgeois contact manifolds arising from open books supporting the standard contact structure on $\mathbb{S}^3$; and similarly for $\mathbb{T}^3$.
Further instances of different open books supporting the same contact structure that induce non-contactomorphic Bourgeois contact manifolds can be found in \cite[Example 1.5]{LMN}; in the same spirit, other examples also come from \Cref{thm:homology_injection_page} proved above.

\begin{question}
Can we find further contactomorphisms of Bourgeois contact manifolds, beyond the inversion of the monodromy from \cite{LMN}? More ambitiously, can we classify the contactomorphism type of all the Bourgeois contact manifolds arising from some fixed contact structure, especially via rigid holomorphic curves invariants? 
\end{question}

In general, \Cref{thm:homology_injection_page} imposes strong constraints on the monodromy. This suggests the following:

\begin{question}
If $BO(\Sigma,\phi)$ is strongly fillable, is $\phi$ (at least smoothly) trivial? 
\end{question}

\bibliographystyle{spmpsci}      

\bibliography{biblio}

\end{document}